\documentclass[11pt]{article}
\usepackage{mathrsfs}
\usepackage{amssymb}
\usepackage{amsmath}
\usepackage{amsbsy}
\usepackage{epsfig}
\usepackage{enumerate}
\usepackage{bm}

\usepackage{color}

\newcommand{\blue}[1]{\textcolor{blue}{#1}}

\usepackage[colorlinks,linkcolor=blue]{hyperref}

\newcommand{\ARXIV}[1]{\href{http://arXiv.org/abs/#1}{\blue{arXiv:#1}}}

\newcommand{\MR}[1]{\href{http://www.ams.org/mathscinet-getitem?mr=#1}{\blue{MR-#1}}}%

\newtheorem{lemma}{Lemma}[section]
\newtheorem{theorem}{Theorem}[section]
\newtheorem{proposition}{Proposition}[section]

\newtheorem{corollary}{Corollary}[section]
\newtheorem{definition}{Definition}[section]
\newtheorem{remark}{Remark}[section]
\newtheorem{example}{Example}[section]

\newbox\TempBox \newbox\TempBoxA

\def\ep{\textsf{E}} 
\def\Sbep{\widehat{\mathbb E}} 
\def\cSbep{\widehat{\mathcal E}} 

\def\Capc{\mathbb V} 
\def\cCapc{\mathcal V} 

\def\underwiggle 1{
\ifmmode\setbox\TempBox=\hbox{$ 1$}\else\setbox\TempBox=\hbox{
1}\fi \setbox\TempBoxA=\hbox to \wd\TempBox{\hss\char'176\hss}
\rlap{\copy\TempBox}\smash{\lower9pt\hbox{\copy\TempBoxA}} }
\renewcommand{\baselinestretch}{1.5}

\begin{document}

\thispagestyle{empty}

\begin{center}
 { \LARGE\bf Lindeberg's central limit theorems for martingale like sequences  under sub-linear expectations$^{\ast}$}
\end{center}

\begin{center} {\sc
\href{https://person.zju.edu.cn/en/stazlx}{\blue{Li-Xin Zhang}}\footnote{Research supported by grants from the NSF of China (Grant No. 11731012),   Ten Thousands Talents Plan of Zhejiang Province (Grant No. 2018R52042), the Fundamental
Research Funds for the Central Universities and the 973 Program
(No. 2015CB352302).
}
}\\
{\sl \small School  of Mathematical Sciences, Zhejiang University, Hangzhou 310027} \\
(Email:stazlx@zju.edu.cn)    \\
\end{center}

\begin{abstract}
The central limit theorem  of martingales is the fundamental tool for studying the convergence of stochastic processes, especially  stochastic integrals and differential equations. In this paper, the central limit theorem and functional central limit theorem are obtained for martingale like random variables  under the sub-linear expectation. As applications, the Lindeberg central limit theorem is obtained for independent but not necessarily identically distributed random variables, and a new proof of the L\'evy characterization of a G-Brownian motion without using stochastic calculus is given. For proving the results,     Rosenthal's inequality and the exponential inequality for the   martingale like random variables are   established.


\end{abstract}

\
\baselineskip 22pt

\renewcommand{\baselinestretch}{1.7}




\section{ Introduction and notations.}\label{sect1}
\setcounter{equation}{0}

Non-additive probabilities and non-additive expectations  are useful tools for studying  uncertainties in statistics, measures of risk,
superhedging in finance and non-linear stochastic calculus, cf. Denis and  Martini (2006), Gilboa (1987),
Marinacci (1999),  Peng (1997, 1999, 2007a, 2007c, 2008a) etc.
Peng  (2007a)   introduced the notion of the sub-linear expectation. Under the sub-linear expectation, Peng (2007a, 2007b, 2007c, 2008a, 2008b, 2009) gave the notions of the G-normal
distributions, G-Brownian motions, G-martingales, independence of random variables, identical distribution of random variables and so on, and developed the weak law of large numbers and central limit theorem for independent and identically distributed (i.i.d.) random variables.
Furthermore, Peng established the stochastic calculus with respect to the G-Brownian motion.  As a result, Peng's   framework of nonlinear expectation gives  a generalization of Kolmogorov's probability theory. Recently, Bayraktar and Munk (2016) proved   an $\alpha$-stable central limit theorem
for independent and identically distributed  random variables.
This paper considers the general  central limit theorem for random variables which are  not necessarily i.i.d.   under the sub-linear expectation. We   establish a central limit theorem and a functional central limit theorem  under the conditional Lindeberg  condition for a kind of martingale-difference like random variables. As   applications,   the central limit theorem  for independent but not necessary identically distributed under the   popular Lindeberg's condition is obtained.    The tool for proving the central limit theorem is a promotion of Peng (2008b)'s and gives also a   new normal approximation method for classical martingale differences instead of the characteristic function. For proving the functional central limit theorem, we  also establish  the Rosenthal's inequalities   for the   martingale like random variables. As the central limit theorem  of classical martingales which is the fundamental tool for studying the convergence of stochastic processes under the  framework of the probability and linear expectation, especially  stochastic integrals and differential equations (cf. Jacod and Shiryaev,2003), the (functional) central limit theorem of martingale-difference like random variables under the sub-linear expectation will provide a way to study the weak convergence of stochastic integrals and difference equations with respect to the G-Brownian motion.

In the rest of this section, we state some notations about   sub-linear expectations. The main results on the central limit theorem and functional central limit theorem are stated in Sections \ref{sectMain} and \ref{sectMain2} with the proofs given the last section. In Section \ref{sectinequality}, we will establish the Rosenthal-type  inequalities and  an exponential   inequality for the maximal sums of the   martingale-difference like random variables.
In Section \ref{sectLevy}, we consider   the L\'evy characterization of a G-Brownian motion in a general sub-linear expectation space. The L\'evy characterization of a $G$-Brownian motion under G-expectation   in a Wiener space is   established by Xu and Zhang (2009, 2010) and extended by Lin (2013)   by the method of  the stochastic calculus.     We will give an elementary proof without using stochastic calculus.  We will find that the functional central limit theorem gives a new way to show  the   L\'evy characterization.

We use the framework and notations of Peng (2008b). Let  $(\Omega,\mathcal F)$
 be a given measurable space  and let $\mathscr{H}$ be a linear space of real functions
defined on $(\Omega,\mathcal F)$ such that if $X_1,\ldots, X_n \in \mathscr{H}$  then $\varphi(X_1,\ldots,X_n)\in \mathscr{H}$ for each
$\varphi\in C_{l,Lip}(\mathbb R_n)$,  where $C_{l,Lip}(\mathbb R_n)$ denotes the linear space of (local Lipschitz)
functions $\varphi$ satisfying
\begin{eqnarray*} & |\varphi(\bm x) - \varphi(\bm y)| \le  C(1 + |\bm x|^m + |\bm y|^m)|\bm x- \bm y|, \;\; \forall \bm x, \bm y \in \mathbb R_n,&\\
& \text {for some }  C > 0, m \in \mathbb  N \text{ depending on } \varphi. &
\end{eqnarray*}
$\mathscr{H}$ is considered as a space of ``random variables''. In this case, we denote $X\in \mathscr{H}$. We also denote   the   space of bounded Lipschitz
functions and the  space of bounded continuous functions on $\mathbb R_n$ by $C_{b,Lip}(\mathbb R_n)$ and  $C_b(\mathbb R_n)$, respectively.

\begin{definition}\label{def1.1} A  sub-linear expectation $\Sbep$ on $\mathscr{H}$  is a function $\Sbep: \mathscr{H}\to \overline{\mathbb R}$ satisfying the following properties: for all $X, Y \in \mathscr H$,
\begin{description}
  \item[\rm (1)]  Monotonicity: If $X \ge  Y$ then $\Sbep [X]\ge \Sbep [Y]$;
\item[\rm (2)] Constant preserving: $\Sbep [c] = c$;
\item[\rm (3)] Sub-additivity: $\Sbep[X+Y]\le \Sbep [X] +\Sbep [Y ]$ whenever $\Sbep [X] +\Sbep [Y ]$ is not of the form $+\infty-\infty$ or $-\infty+\infty$;
\item[\rm (4)] Positive homogeneity: $\Sbep [\lambda X] = \lambda \Sbep  [X]$, $\lambda\ge 0$.
 \end{description}
 Here $\overline{\mathbb R}=[-\infty, \infty]$. The triple $(\Omega, \mathscr{H}, \Sbep)$ is called a sub-linear expectation space. Give a sub-linear expectation $\Sbep $, let us denote the conjugate expectation $\cSbep$ of $\Sbep$ by
$  \cSbep[X]:=-\Sbep[-X]$,  $ \forall X\in \mathscr{H}$.
\end{definition}
A sub-linear expectation $\Sbep$ is countably sub-additive, if
$$ \Sbep[\sum_{i=1}^{\infty}X_i]\le \sum_{i=1}^{\infty} \Sbep[X_i], \;\; \text{for all random variables } X_i\ge 0. $$
If $X$ is not in $\mathscr{H}$, we define its sub-linear expectation by $\Sbep^{\ast}[X]=\inf\{\Sbep[Y]:   X\le Y\in \mathscr{H}\}$.  When there is no ambiguity, we also denote it by  $\Sbep$.
From the definition, it is easily shown that    $\cSbep[X]\le \Sbep[X]$, $\Sbep[X+c]= \Sbep[X]+c$ and $\Sbep[X-Y]\ge \Sbep[X]-\Sbep[Y]$ for all
$X, Y\in \mathscr{H}$ with $\Sbep[Y]$ being finite. Further, if $\Sbep[|X|]$ is finite, then $\cSbep[X]$ and $\Sbep[X]$ are both finite.

\begin{definition} ({\em Peng (2007a, 2008b)})

\begin{description}
  \item[ \rm (i)] ({\em Identical distribution}) Let $\bm X_1$ and $\bm X_2$ be two $n$-dimensional random vectors defined,
respectively, in sub-linear expectation spaces $(\Omega_1, \mathscr{H}_1, \Sbep_1)$
  and $(\Omega_2, \mathscr{H}_2, \Sbep_2)$. They are called identically distributed, denoted by $\bm X_1\overset{d}= \bm X_2$,  if
$$ \Sbep_1[\varphi(\bm X_1)]=\Sbep_2[\varphi(\bm X_2)], \;\; \forall \varphi\in C_{l,Lip}(\mathbb R_n), $$
whenever the sub-expectations are finite. A sequence $\{X_n;n\ge 1\}$ of random variables is said to be identically distributed if $X_i\overset{d}= X_1$ for each $i\ge 1$.
\item[\rm (ii)] ({\em Independence})   In a sub-linear expectation space  $(\Omega, \mathscr{H}, \Sbep)$, a random vector $\bm Y =
(Y_1, \ldots, Y_n)$, $Y_i \in \mathscr{H}$ is said to be independent to another random vector $\bm X =
(X_1, \ldots, X_m)$ , $X_i \in \mathscr{H}$ under $\Sbep$,  if for each test function $\varphi\in C_{l,Lip}(\mathbb R_m \times \mathbb R_n)$
we have
$ \Sbep [\varphi(\bm X, \bm Y )] = \Sbep \big[\Sbep[\varphi(\bm x, \bm Y )]\big|_{\bm x=\bm X}\big],$
whenever $\overline{\varphi}(\bm x):=\Sbep\left[|\varphi(\bm x, \bm Y )|\right]<\infty$ for all $\bm x$ and
 $\Sbep\left[|\overline{\varphi}(\bm X)|\right]<\infty$.

 Random variables $X_1,\ldots, X_n$ are said to be independent if for each $2\le k\le n$, $X_k$ is independent to $(X_1,\ldots, X_{k-1})$. A sequence of random variables is said to be independent if for each $n$, $X_1,\ldots, X_n$ are independent.
\end{description}
\end{definition}

Next, we introduce the capacities corresponding to the sub-linear expectation.
  We denote the pair $(\Capc,\cCapc)$ of capacities on  $(\Omega, \mathscr{H}, \Sbep)$ by setting
$$ \Capc(A):=\inf\{\Sbep[\xi]: I_A\le \xi, \xi\in\mathscr{H}\}, \;\; \cCapc(A):= 1-\Capc(A^c),\;\; \forall A\in \mathcal F, $$
where $A^c$  is the complement set of $A$.
Then
it is obvious that $\Capc$ is sub-additive, i.e. $\Capc(A\bigcup B)\le \Capc(A)+\Capc(B)$. But $\cCapc$ and $\cSbep$ are not. However, we have
$$
  \cCapc(A\bigcup B)\le \cCapc(A)+\Capc(B) \;\;\text{ and }\;\; \cSbep[X+Y]\le \cSbep[X]+\Sbep[Y]
$$
due to the fact that $\Capc(A^c\bigcap B^c)=\Capc(A^c\backslash B)\ge \Capc(A^c)-\Capc(B)$ and $\Sbep[-X-Y]\ge \Sbep[-X]-\Sbep[Y]$.

The  Choquet integrals/expecations of  $(C_{\Capc},C_{\cCapc})$  are defined by
$$ C_V[X]=\int_0^{\infty} V(X\ge t)dt +\int_{-\infty}^0\left[V(X\ge t)-1\right]dt $$
with $V$ being replaced by $\Capc$ and $\cCapc$, respectively.

Finally, we recall the notations of G-normal distribution and  G-Brownian motion which are introduced by Peng (2008b, 2010).

\begin{definition} ({\em G-normal random variable})
For $0\le \underline{\sigma}^2\le \overline{\sigma}^2<\infty$, a random variable $\xi$ in a sub-linear expectation space $(\widetilde{\Omega}, \widetilde{\mathscr H}, \widetilde{\mathbb E})$   is called a normal $N\big(0, [\underline{\sigma}^2, \overline{\sigma}^2]\big)$ distributed  random variable   (written as $\xi\sim N\big(0, [\underline{\sigma}^2, \overline{\sigma}^2]\big)$  under $\widetilde{\mathbb E}$), if for any  $\varphi\in C_{l,Lip}(\mathbb R)$, the function $u(x,t)=\widetilde{\mathbb E}\left[\varphi\left(x+\sqrt{t} \xi\right)\right]$ ($x\in \mathbb R, t\ge 0$) is the unique viscosity solution of  the following heat equation:
      $$ \partial_t u -G\left( \partial_{xx}^2 u\right) =0, \;\; u(0,x)=\varphi(x), $$
where $G(\alpha)=\frac{1}{2}(\overline{\sigma}^2 \alpha^+ - \underline{\sigma}^2 \alpha^-)$.
\end{definition}
That $\xi$ is a normal distributed random variable is equivalent to that,   if $\xi^{\prime}$ is an independent copy of $\xi$, then
$$ \widetilde{\mathbb E}\left[\varphi(\alpha \xi+\beta \xi^{\prime})\right]
=\widetilde{\mathbb E}\left[\varphi\big(\sqrt{\alpha^2+\beta^2}X\big)\right], \;\;
\forall \varphi\in C_{l,Lip}(\mathbb R) \text{ and } \forall \alpha,\beta\ge 0, $$
(cf. Definition II.1.4 and Example II.1.13 of Peng (2010)).

\begin{definition}\label{DefG-B} ({\em $G$-Brownian motion})  A random process $(W_t)_{t\ge 0}$ in the sub-linear expectation space $(\widetilde{\Omega}, \widetilde{\mathscr H}, \widetilde{\mathbb E})$ is called a $G$-Brownian motion (cf. Definition III.1.2 of Peng (2010)) if
\begin{description}
  \item[\rm (i)] $W_0=0$;
  \item[\rm (ii)]  For each $0\le t_1\le \cdots\le t_d\le t\le s$,
\begin{align}
&\widetilde{\mathbb E}\left[\varphi\big(W_{t_1},\ldots, W_{t_d}, W_s-W_t\big)\right]\nonumber
\\
= &
  \widetilde{\mathbb E}\left[\widetilde{\mathbb E}\left[\varphi\big(x_1,\ldots, x_d, \sqrt{t-s})\xi\big)\right]\big|_{x_1=W_{t_1},\ldots, x_d=W_{t_d}}\right]
  \label{eqBrown} \\
  & \;\; \forall \varphi\in C_{l,Lip}(\mathbb R_{d+1}), \nonumber
  \end{align}
  where $\xi\sim N(0,[\underline{\sigma}^2,\overline{\sigma}^2])$.
\end{description}
\end{definition}
In some papers, for example, Xu and Zhang (2009, 2010),  the test functions $\varphi$ are only required to be elements  in  $C_{b,Lip}(\mathbb R_{d+1})$. It can be shown that if $\widetilde{\mathbb E}[|W_t|^p]<\infty$ for all $p>0$ and $t$, then that (\ref{eqBrown}) holds for all $\varphi\in C_{b,Lip}(\mathbb R_{d+1})$ is equivalent to that it holds for all $\varphi\in C_{l,Lip}(\mathbb R_{d+1})$. Further, if the sub-linear expectation $\widetilde{\mathbb E}$ is countably sub-additive, then this two kinds of definitions are equivalent because, if    $X$ is a random variable in $(\Omega,\mathscr{H}, \Sbep)$ such that
\begin{equation}\label{boundedID} \Sbep[\varphi(X)]=\widetilde{\mathbb E} [\varphi(\xi)], \;\; \forall  \varphi\in C_{b,Lip}(\mathbb R),
\end{equation}
then $\Sbep[|X|^p]<\infty$ for all $p>0$.
 In fact,
if $\xi\sim N\big(0, [\underline{\sigma}^2, \overline{\sigma}^2]\big)$ under $\widetilde{\mathbb E}$, then (cf. Peng(2010, page 22))
$$ \widetilde{\mathbb E}[|\xi|^p]=\overline{\sigma}^p\int_{-\infty}^{\infty} |x|^p \frac{1}{\sqrt{2\pi}}e^{-x^2/2} dx=c_p\overline{\sigma}^p, \;\; \forall p\ge 1. $$
Now,
for any $z>0$, one can choose a function $\varphi\in C_{b,Lip}(\mathbb R)$ such that $I\{x> z\}\le \varphi(x)\le I\{x> z-\epsilon\}$. From (\ref{boundedID}), it follows that
$$ \Capc(|X|>z)\le \Sbep[\varphi(X)]=\widetilde{\mathbb E} [\varphi(\xi)]\le \widetilde{\Capc}(|\xi|>z-\epsilon). $$
 Hence
$$ \Capc(|X|>z)\le   \widetilde{\Capc}(|\xi|\ge z/2)\le \frac{2^{2p}\widetilde{\mathbb E}[|\xi|^{2p}]}{z^{2p}}=\frac{\overline{\sigma}^{2p}c_{2p}}{z^{2p}}. $$
It follows that
$$ C_{\Capc}(|X|^p)=\int_0^{\infty} \Capc(|X|^p>z)dz\le 1+\int_1^{\infty}\frac{\overline{\sigma}^{2p}c_{2p}}{z^2}dz\le 1+\overline{\sigma}^{2p}c_{2p}
<\infty, \;\; \forall p\ge 2. $$
So, if $\Sbep$ is countably sub-additive or $\Sbep[|X|^p]=\lim_{c\to \infty} \Sbep[(|X|\wedge c)^p]$,
then $\Sbep[|X|^p]\le C_{\Capc}(|X|^p)<\infty$ for all $p>0$ by Lemma 3.9 of Zhang (2016).

Let $C_{[0,1]}$ be a function space of continuous real functions   on $[0,1]$ equipped with the supremum norm $\|x\|=\sup_{0\le t\le 1}|x(t)|$ and $C_b\big(C_{[0,1]}\big)$ is the set of bounded continuous  functions $h(x):C_{[0,1]}\to \mathbb R$. As showed in  Peng (2006, 2008a, 2010) and Denis, Hu, and Peng (2011),    there is a sub-linear expectation space $\big(\widetilde{\Omega}, \widetilde{\mathscr{H}},\widetilde{\mathbb E}\big)$ with
$\widetilde{\Omega}= C_{[0,1]}$ and $C_b\big(C_{[0,1]}\big)\subset \widetilde{\mathscr{H}}$ such that $(\widetilde{\mathscr{H}}, \widetilde{\mathbb E}[\|\cdot\|])$ is a Banach space, and
the canonical process $W(t)(\omega) = \omega_t  (\omega\in \widetilde{\Omega})$ is a G-Brownian motion. In the sequel of this paper, the G-normal random variables and G-Brownian motions are considered in $(\widetilde{\Omega}, \widetilde{\mathscr{H}}, \widetilde{\mathbb E})$.

\section{Lindeberg's central limit theorem for independent random variables.}\label{sectMain}
\setcounter{equation}{0}

We write $\eta_n\overset{\Capc}\to \eta$ if $\Capc\left(|\eta_n-\eta|\ge \epsilon\right)\to 0$ for any $\epsilon>0$, and write $\eta_n\overset{d}\to \eta$ if
$\Sbep\left[\varphi(\eta_n)\right]\to \Sbep\left[\varphi(\eta)\right]$ holds for all bounded and continuous functions $\varphi$.
 In this section , we consider the independent random variables $\{X_{n,k};k=1,\ldots, k_n\}$.
Denote    $\overline{\sigma}_{n,k}^2=\Sbep[X_{n,k}^2]$,  $\underline{\sigma}_{n,k}^2=\cSbep[X_{n,k}^2]$, $B_n^2=\sum_{k=1}^{k_n} \overline{\sigma}_{n,k}^2$ . We have the following Lindeberg's central limit theorem.
\begin{theorem}\label{thCLT}    Suppose that the Lindeberg condition is satisfied:
 \begin{align}\label{eqLindeberg}
 \frac{1}{B_n^2}\sum_{k=1}^{k_n}\Sbep\left[ \left( X_{n,k}^2-\epsilon B_n^2 \right)^+\right]\to 0\;\; \forall \epsilon>0,
 \end{align}
 and further, there is a constant $r\in[0,1]$ such that
  \begin{align}\label{eqCLTCond.2}
 \frac{\sum_{k=1}^{k_n} \left| r\overline{\sigma}_{n,k}^2 - \underline{\sigma}_{n,k}^2\right| }{B_n^2}\to 0, \;\; \text{also, }
 \end{align}
  \begin{align}\label{eqCLTCond.3}
 \frac{\sum_{k=1}^{k_n}\left\{ |\Sbep[X_{n,k}]|+|\cSbep[X_{n,k}]|\right\}  }{B_n}\to 0.
 \end{align}
 Then for any  bounded continuous function $\varphi$,
  \begin{equation} \label{eqCLT} \lim_{n\to \infty}\Sbep\left[\varphi\left(\frac{\sum_{k=1}^{k_n}X_{n,k}}{B_n}\right)\right]=\widetilde{\mathbb E}[\varphi(\xi)],
  \end{equation}
   where   $\xi\sim N(0,[r, 1])$  under $\widetilde{\mathbb E}$.
 \end{theorem}

Theorem \ref{thCLT} will be a directly corollary of our   Theorem  \ref{thCLTM} on the central limit theorem for martingale like sequence. The central limit theorem for independent and identically distributed random variables under the sub-linear expectation was obtained by Peng (2008b). Li and Shi (2010) generalized Peng's result to a central limit theorem for   independent random variables $\{X_n;n\ge 1 \}$ satisfying $\Sbep[X_i]=\cSbep[X_i]=0$, $\Sbep[|X_i|^3]\le M<\infty$, $i=1,2, \ldots $, and
$$ \frac{1}{n}\sum_{i=1}^n \left|\Sbep[X_i^2]-\overline{\sigma}^2\right|\to 0, \;\;  \frac{1}{n}\sum_{i=1}^n \left|\cSbep[X_i^2]-\underline{\sigma}^2\right|\to 0. $$
It is easily seen that the array $\{\frac{1}{\sqrt{n}} X_k; k=1,\ldots, n\}$ satisfies the conditions (\ref{eqCLTCond.2})  with   $r=\underline{\sigma}^2/\overline{\sigma}^2$,   (\ref{eqCLTCond.3}) and
(\ref{eqLindeberg}).

  When $\Sbep$ is a classical linear expectation,  (\ref{eqCLTCond.2}) is automatically satisfied with $r=1$.  It is easily seen that (\ref{eqCLTCond.2}) implies
\begin{equation}\label{eqCLTCond.2ad}  \frac{\;\;\sum_{k=1}^{k_n}  \underline{\sigma}_{n,k}^2\;\;}{\sum_{k=1}^{k_n}  \overline{\sigma}_{n,k}^2} \to r.
\end{equation}
One may conjecture that (\ref{eqCLTCond.2}) can be weakened to (\ref{eqCLTCond.2ad}). The following example tells us that it is not the truth.

\begin{example} Let $0<\tau_1,\tau_2<1$,  and $\{X_{n,k}; k=1,\ldots, 2n\}$ be a sequence of independent normal random variables   such that
$$ X_{n,k} \overset{d}\sim N\big(0,[\tau_1,1]), k=1,\ldots, n \text{ and }  X_{n,k} \overset{d}\sim    N\big(0,[\tau_2,1]), k=n+1,\ldots, 2n. $$
It is easily seen that $\{X_{n,k}; k=1,\ldots, 2n\}$  satisfies the conditions (\ref{eqLindeberg}), (\ref{eqCLTCond.3}) and (\ref{eqCLTCond.2ad}) with
 $r=(\tau_1+\tau_2)/2$, and $B_n^2= 2 n$. It is obvious that
$$ \frac{\sum_{k=1}^{2n} X_{n,k}}{\sqrt{n}}=\frac{\sum_{k=1}^{n} X_{n,k}}{\sqrt{n}}+\frac{\sum_{k=n+1}^{2n} X_{n,k}}{\sqrt{n}}\overset{d}\sim \xi +  \eta, $$
where $\xi,\eta$ are independent normal random variables with $\xi\overset{d}\sim N(0,[\tau_1,1])$, $\eta\overset{d}\sim N(0,[\tau_2,1])$.  Song (2015) showed that
$\xi +  \eta$ is not   $G$-normal distributed if $\tau_1\ne \tau_2$, and hence (\ref{eqCLT}) fails.
\end{example}

\section{Central limit theorem for martingale like sequence.}\label{sectMain2}
\setcounter{equation}{0}

 In this section, we consider a general martingale. First, we  recall the definition of the conditional expectation under the sub-linear expectation. Let $(\Omega, \mathscr{H}, \Sbep)$ be a sub-linear expectation space.  We write $X\le Y$ in $L_p$ if $\Sbep[((X-Y)^+)^p]=0$, $X= Y$ in $L_p$ if both $X\le Y$ and $Y\le X$ holds in $L_p$.

  Let   $\mathscr{H}_{n,0}\subset \cdots\subset
 \mathscr{H}_{n,k_n}$ be subspaces of $\mathscr{H}$ such that
  \begin{description}
    \item[\rm (1)]  any constant   $c\in \mathscr{H}_{n,k}$ and,
    \item[\rm (2)] if $X_1,\ldots,X_d\in \mathscr{H}_{n,k}$, then $\varphi(X_1,\ldots,X_d)\in \mathscr{H}_{n,k}$ for any $\varphi\in C_{l,lip}(\mathbb R_d)$, $k=0,\cdots, k_n$.
  \end{description}
  Denote $\mathscr{L}(\mathscr{H})=\{X:\Sbep[|X|]<\infty, X\in \mathscr{H}\}$.
 We consider a system of operators in $\mathscr{L}(\mathscr{H})$,
 $$\Sbep_{n,k}: \mathscr{L}(\mathscr{H})\to \mathscr{L}(\mathscr{H}_{n,k}) $$
 and denote $\Sbep[X|\mathscr{H}_{n,k}]=\Sbep_{n,k}[X]$, $\cSbep[X|\mathscr{H}_{n,k}]=-\Sbep_{n,k}[-X]$.  $\Sbep[X|\mathscr{H}_{n,k}]$ is called the conditional sub-linear expectation of $X$ given $\mathscr{H}_{n,k}$, $\Sbep_{n,k}$ is called the conditional expectation operator.
Suppose that the operators $\Sbep_{n,k}$ satisfy  the following properties:  for all $X, Y \in \mathscr{L}({\mathscr H})$,
\begin{description}
  \item[\rm (a)]   $ \Sbep_{n,k} [ X+Y]=X+\Sbep_{n,k}[Y]$ in $L_1$ if $X\in \mathscr{H}_{n,k}$, and $ \Sbep_{n,k} [ XY]=X^+\Sbep_{n,k}[Y]+X^-\Sbep_{n,k}[-Y]$ in $L_1$  if
$X\in \mathscr{H}_{n,k}$ and $XY\in \mathscr{L}({\mathscr H})$;
\item[\rm (b)]   $\Sbep\left[\Sbep_{n,k} [ X]\right]=\Sbep[X]$.
 \end{description}
It is easily seen that (a) implies that  $\Sbep_{n,k} [c] = c$ in $L_1$ and   $\Sbep_{n,k} [\lambda X] = \lambda \Sbep_{n,k}  [X]$ in $L_1$ if $\lambda\ge 0$.
The definition of the conditional sub-linear expectation can be found in  Peng (2010), Xu and Zhang (2009, 2010)  with the operators satisfying (a), (b) and,
  $\Sbep_{n,k}[X]\le   \Sbep_{n,k}[Y]$ if $X\le Y$,
 $\Sbep_{n,k}[X]-\Sbep_{n,k}[Y]\le \Sbep_{n,k}[X-Y]$,
 $\Sbep_{n,k}\left[\left[\Sbep_{n,l} [ X]\right]\right]=\Sbep_{n,l\wedge k} [ X]$.
It can be showed that these properties can be implied by (a) and (b) (c.f. Lemma \ref{lemma4.0.2}).

Now, we assume that $\{Z_{n,k}; k=1,\ldots, k_n\}$ is an array of random variables such that $Z_{n,k}\in \mathscr{H}_{n,k}$ and $\Sbep[Z_{n,k}^2]<\infty$, $k=1,\ldots, k_n$. The following is the central limit theorem.
 \begin{theorem} \label{thCLTM}   Suppose that the operators $\Sbep_{n,k}$ satisfy (a) and (b).   Assume that   the following Lindeberg condition is satisfied:
 \begin{equation}\label{eqLindebergM}
  \sum_{k=1}^{k_n}\Sbep\left[ \left( Z_{n,k}^2-\epsilon \right)^+|\mathscr{H}_{n,k-1}\right]\overset{\Capc}\to 0\;\; \forall \epsilon>0,
 \end{equation}
 and further, there are   constants $\rho\ge 0$ and  $r\in[0,1]$ such that
 \begin{equation}\label{eqCLTCondM.2}
 \sum_{k=1}^{k_n}  \Sbep[Z_{n,k}^2|\mathscr{H}_{n,k-1}]    \overset{\Capc}\to \rho,
 \end{equation}
  \begin{equation}\label{eqCLTCondM.3}
 \sum_{k=1}^{k_n} \left| r\Sbep[Z_{n,k}^2|\mathscr{H}_{n,k-1}] - \cSbep[Z_{n,k}^2|\mathscr{H}_{n,k-1}]\right|  \overset{\Capc}\to 0,
 \end{equation}
  \begin{equation}\label{eqCLTCondM.4}
  \sum_{k=1}^{k_n}\left\{ |\Sbep[Z_{n,k} |\mathscr{H}_{n,k-1}]|+|\cSbep[Z_{n,k} |\mathscr{H}_{n,k-1}]|\right\}    \overset{\Capc}\to 0.
 \end{equation}
Then for any  bounded   continuous function $\varphi$,
  \begin{equation} \label{eqCLTM} \lim_{n\to \infty}\Sbep\left[\varphi\left(\sum_{k=1}^{k_n}Z_{n,k}\right)\right]=\widetilde{\mathbb E}[\varphi(\sqrt{\rho}\xi)],
  \end{equation}
  i.e.,
   $   \sum_{k=1}^{k_n}Z_{n,k}\overset{d}\to \sqrt{ \rho} \xi,$
     where   $\xi\sim N(0,[r, 1])$  under $\widetilde{\mathbb E}$.
 \end{theorem}

\begin{remark} When $ \Sbep[Z_{n,k}|\mathscr{H}_{n,k-1}]=0$ and $\cSbep[Z_{n,k}|\mathscr{H}_{n,k-1}]=0$, then $\{Z_{n,k}; k=1,\ldots, k_n\}$ is an array of symmetric martingale differences (cf. Xu and Zhang (2009)). If $\Sbep[\cdot]=\ep_P[\cdot]$ is a classical linear expectation, then  (\ref{eqCLTCondM.3}) is satisfied with $r=1$, and the conclusion coincides with Corollary 3.1 of Hall and Heyde (1980).
\end{remark}

The following is a direct corollary of Theorem \ref{thCLTM}.

\begin{corollary}\label{Cor1}
Let $\{\eta_n\}$ be a sequence of independent   random variables on $(\Omega, \mathscr{H}, \Sbep)$ with $\Sbep[\eta_n]=\cSbep[\eta_n]=0$, $\Sbep[\eta_n^2]=:\overline{\sigma}^2_n\to \overline{\sigma}^2$,  $\cSbep[\eta_n^2]:=\underline{\sigma}_n^2\to \underline{\sigma}^2$ and $\sup_n\Sbep[(\eta_n^2-c)^+ ]\to 0 $ as $c\to \infty$.
Suppose that $\{a_{n,i}; i=1,\ldots, k_n\}$ is an array of real random variables in $\mathscr{H}$   with $a_{n,i}$ being a function of $\eta_1,\ldots,\eta_{i-1}$,
$$\max_i |a_{n,i}|\overset{\Capc}\to 0  \text{ and } \sum_{i=1}^{k_n} a_{n,i}^2\overset{\Capc}\to \rho, $$
where $\rho\ge 0$ is a constant.  Then
\begin{equation}\label{eqCor1.1} \lim_{n\to \infty} \Sbep\left[\varphi\Big(\sum_{i=1}^{k_n}a_{n,i}\eta_i\Big)\right]=\widetilde{\mathbb E}[\varphi(\xi)],
  \end{equation}
for any bounded   continuous function $\varphi$,  where   $\xi\sim N(0,[\rho\underline{\sigma}^2, \rho\overline{\sigma}^2])$  under $\widetilde{\mathbb E}$.
\end{corollary}

The following corollary is a central limit theorem for moving average processes which include the ARMA model.
\begin{corollary}Let $\{\eta_n\}$ be a sequence of independent and identically distributed   random variables in $(\Omega, \mathscr{H}, \Sbep)$ with $\Sbep[\eta_1]=\cSbep[\eta_1]=0$, $\Sbep[\eta_1^2]=\overline{\sigma}^2$ and $\cSbep[\eta_1^2]=\underline{\sigma}^2$, $\{a_n;n\ge 0\}$ be a sequence of real numbers
with $\sum_{n=0}^{\infty}|a_n|<\infty$. Let
$X_k=\sum_{i=0}^{\infty}a_i\eta_{i+k}. $
Then
\begin{equation}\label{eqCor2.1}
\frac{1}{\sqrt{n}}\sum_{k=1}^n X_k \overset{d}\to N\big(0,[a^2\underline{\sigma}^2, a^2\overline{\sigma}^2]\big),
\end{equation}
where $a=\sum_{j=0}^{\infty}a_j$.
\end{corollary}
{\bf Proof.} Let $a_n=0$ if $n<0$. Then $X_k= \sum_{i=1}^{\infty} a_{i-k}\eta_i$ and
$$ \frac{1}{\sqrt{n}}\sum_{k=1}^n X_k =\sum_{i=1}^{\infty} \left(\frac{\sum_{k=1}^na_{i-k}}{\sqrt{n}}\right)\eta_i. $$
 Let $a_{n,i}=\frac{\sum_{k=1}^na_{i-k}}{\sqrt{n}}$. Then $\max_i|a_{n,i}|\le n^{-1/2}\sum_{i=-\infty}^{\infty}|a_i|\to 0$ and
 $ \sum_{i=1}^{\infty} a_{n,i}^2\to a^2. $
The result follows from Corollary \ref{Cor1}. $\Box$

\bigskip
Finally, we give the functional central limit theorems.

Let $D_{[0,1]}$ be the space of right continuous functions having finite left limits which is endowed with the Skorohod topology, $\tau_n(t)$ be a non-decreasing function in $D_{[0,1]}$ which takes integer values with $\tau_n(0)=0$, $\tau_n(1)=k_n$. Define  $S_{n,i}=\sum_{k=1}^i Z_{n,k}$,
\begin{equation}\label{eqthFCLTM.1} W_n(t)= S_{n, \tau_n(t)}.
\end{equation}

 \begin{theorem}\label{thFCLTM}   Suppose that the operators $\Sbep_{n,k}$ satisfy (a) and (b).   Assume that   the conditions (\ref{eqLindebergM}), (\ref{eqCLTCondM.3}) and (\ref{eqCLTCondM.4}) in Theorem \ref{thCLTM} are satisfied. Further,   there is  a continuous non-decreasing non-random function   $\rho(t)$  such that
 \begin{equation}\label{eqFCLTCondM.2}
 \sum_{k\le \tau_n(t)}  \Sbep[Z_{n,k}^2|\mathscr{H}_{n,k-1}]    \overset{\Capc}\to \rho(t), \;\; t\in[0,1].
 \end{equation}
 Then  for any $0=t_0<\ldots< t_d\le 1$,
 \begin{equation} \label{eqfinitedimension} \Big( W_n(t_1),\cdots, W_n(t_d)\Big)\overset{d}\to \Big( W(\rho(t_1)),\cdots, W(\rho(t_d))\Big),
 \end{equation}
 and for any  bounded continuous function $\varphi:D_{[0,1]}\to \mathbb R$,
  \begin{equation} \label{eqFCLTM} \lim_{n\to \infty}\Sbep\left[\varphi\left(W_n\right)\right]=\widetilde{\mathbb E}[\varphi(W\circ\rho)],
  \end{equation}
   where  $W$ is $G$-Brownian motion on $[0,1]$ with  $W(1) \sim N(0,[r, 1])$  under $\widetilde{\mathbb E}$, and $W\circ\rho(t)=W(\rho(t))$.
 \end{theorem}

Because the proofs of Theorems \ref{thCLTM} and \ref{thFCLTM} are a little long and need some preparation, we will give them in the last section.

  \section{Moment inequalities and exponential   inequalities.}\label{sectinequality}
  \setcounter{equation}{0}

To prove the central limit theorems and functional central limit theorems, we need some inequalities on the sums of martingale-difference like random variables  as basic tools.
Before we give the inequalities, we state some properties of the sub-linear expectations $\Sbep$ and $\Sbep_{n,k}$.
  The first is     H\"older's inequality which is Proposition 16 of Denis,  Hu,  and Peng (2011).
 \begin{lemma} ({\em H\"older's inequality}) Let $p,q>1$ be two real numbers satisfying $\frac{1}p+\frac{1}{q}=1$. Then, for two random variables  $X,Y$  in $(\Omega, \mathscr{H}, \Sbep)$ we have
 $$ \Sbep[|XY|]\le \left(\Sbep[|X|^p]\right)^{\frac{1}{p}}  \left(\Sbep[|Y|^q]\right)^{\frac{1}{q}}$$
whenever $\Sbep[|X|^p]<\infty$, $\Sbep[|Y|^q]<\infty$.
 \end{lemma}

The next two lemmas are on the properties of the sub-linear expectation, the capacity and the operators $\Sbep_{n,k}$. The  proofs will be given   in Appendix \ref{appendixA}.  We write  $X\le Y$ in capacity $\Capc$ if $\Capc\left(X-Y\ge \epsilon\right)=0$ for all $\epsilon>0$,
 and $X= Y$ in capacity $\Capc$ if   both $X\le Y$ and $Y\le X$ holds in $\Capc$
 \begin{lemma}\label{lemma4.0.1} We have
 \begin{description}
   \item[\rm (1)]  if $X\le Y$ in $L_p$, then  $X\le Y$ in   $\Capc$;
   \item[\rm (2)]  if $X\le Y$ in   $\Capc$ and $\Sbep[((X-Y)^+)^p]<\infty$, then $X\le Y$ in $L_q$ for $0<q<p$;
   \item[\rm (3)]  if $X\le Y$ in   $\Capc$, $f(x)$ is non-decreasing continuous function  and $\Capc(|Y|\ge M)\to 0$ as $M\to \infty$, then $f(X)\le f(Y)$ in   $\Capc$;
   \item[\rm (4)] if $p\ge 1$, $X,Y\ge 0$ in $L_p$, $X\le Y$ in $L_p$, then $\Sbep[X^p]\le \Sbep[Y^p]$;
   \item[\rm (5)] if $\Sbep$ is countably additive, then $X\le Y$ in $\Capc$ is equivalent to  $X\le Y$ in $L_p$ for any $p>0$.
 \end{description}
 \end{lemma}

\begin{lemma}\label{lemma4.0.2} Suppose that the operators $\Sbep_{n,k}$ satisfy (a) and (b). For $X,Y\in\mathscr{L}(\mathscr{H})$, we have
\begin{description}
  \item[\rm (c)] if $X\le Y$ in $L_1$, then $\Sbep_{n,k}[X]\le   \Sbep_{n,k}[Y]$ in $L_1$;
 \item[\rm (d)]  $\Sbep_{n,k}[X]-\Sbep_{n,k}[Y]\le \Sbep_{n,k}[X-Y]\le \Sbep_{n,k}[|X-Y|]$ in $L_1$;
\item[\rm (e)]   $\Sbep_{n,k}\left[\left[\Sbep_{n,l} [ X]\right]\right]=\Sbep_{n,l\wedge k} [ X]$ in $L_1$;
\item[\rm (f)] if $|X|\le M$ in $L_p$ for all $p\ge  1$, then $ \big|\Sbep_{n,k}[X]\big| \le M$ in $L_p$ for all $p\ge  1$.
 \end{description}
\end{lemma}

\bigskip
For the martingale-difference like random variables, we have the following  theorem on the  Rosenthal-type inequalities.

\begin{theorem} \label{lemRosenIeq}   Set $S_0=0$, $S_k=\sum_{i=1}^k  Z_{n,i}$. Suppose that $\{Z_{n,i}\}$ are a set of bounded random variables. Then,
\begin{align}\label{eqlemRosenIeq.1}
\Sbep\Big[\Big(\max_{k\le {k_n}} (S_{k_n}-S_k)\Big)^2\Big]\le      \Sbep\Big[\sum_{k=1}^{k_n}   \Sbep[Z_{n,k}^2|\mathscr{H}_{n,k-1}]\Big]
\end{align}
when $\Sbep[Z_{n,k}|\mathscr{H}_{n,k-1}]\le 0$, $k=1,\ldots, k_n$, and in general,
\begin{align}\label{eqlemRosenIeq.2}
\Sbep\Big[\max_{k\le {k_n}} |S_k|^2\Big]\le &  256\left\{\Sbep\Big[\sum_{k=1}^{k_n}   \Sbep[Z_{n,k}^2|\mathscr{H}_{n,k-1}]\Big]\right.\nonumber\\
&   \left.+\Sbep\left[\Big\{\sum_{k=1}^{k_n}  \Big(\big( \Sbep[Z_{n,k} |\mathscr{H}_{n,k-1}]\big)^++\big( \cSbep[Z_{n,k} |\mathscr{H}_{n,k-1}]\big)^-\Big)\Big\}^2\right]\right\}.
\end{align}
Moreover,   for $p\ge 2$ there is a constant $C_p$ such that
\begin{align}\label{eqlemRosenIeq.4}
\Sbep\Big[\max_{k\le {k_n}} |S_k|^p\Big]\le &  C_p\left\{\Sbep\left[\sum_{k=1}^{k_n} \Sbep[|Z_{n,k}|^p|\mathscr{H}_{n,k-1}]\right]+\Sbep\left[\Big(\sum_{k=1}^{k_n}   \Sbep[Z_{n,k}^2|\mathscr{H}_{n,k}]\Big)^{p/2}\right]\right.\nonumber\\
&  \qquad \left.+\Sbep\left[\Big\{\sum_{k=1}^{k_n}  \Big(\big( \Sbep[Z_{n,k} |\mathscr{H}_{n,k}]\big)^++\big( \cSbep[Z_{n,k} |\mathscr{H}_{n,k}]\big)^-\Big)\Big\}^p\right]\right\}.
\end{align}
\end{theorem}
{\bf Proof. } Let $Q_k=\max\{  Z_{n,k},  Z_{n,k}+  Z_{n,k-1}, \ldots,  Z_{n,k}+\cdots  Z_{n,1}\}$,   $M_k=\max_{i\le k}|S_i|$.  Then, $Q_k=  Z_{n,k}+Q_{k-1}^+$, $Q_k^2= Z_{n,k}^2+2 Z_{n,k} Q_{k-1}^++(Q_{k-1}^+)^2$, $|Q_k|\le 2 M_{k_n}$. It follows that
\begin{align*}
 & \left(\max_{k\le k_n}(S_{k_n}-S_k)\right)^2=   (Q_{k_n}^+)^2 \le    \sum_{k=1}^{k_n} Z_{n,k}^2+2 \sum_{k=1}^{k_n} Z_{n,k} Q_{k-1}^+ \\
 \le &\sum_{k=1}^{k_n} \Sbep[Z_{n,k}^2|\mathscr{H}_{n,k-1}]+\sum_{k=1}^{k_n}  \big(Z_{n,k}^2-\Sbep[Z_{n,k}^2|\mathscr{H}_{n,k-1}]) \\
& +2\sum_{k=1}^{k_n} \Sbep[Z_{n,k}|\mathscr{H}_{n,k-1}] Q_{k-1}^+ +2\sum_{k=1}^{k_n}\big(Z_{n,k}- \Sbep[Z_{n,k}|\mathscr{H}_{n,k-1}]\big) Q_{k-1}^+\\
  \le &\sum_{k=1}^{k_n} \Sbep[Z_{n,k}^2|\mathscr{H}_{n,k-1}]+4\sum_{k=1}^{k_n} (\Sbep[Z_{n,k}|\mathscr{H}_{n,k-1}])^+ M_{k_n} \\
  &+\sum_{k=1}^{k_n}  \big(Z_{n,k}^2-\Sbep[Z_{n,k}^2|\mathscr{H}_{n,k-1}])
   +2\sum_{k=1}^{k_n}\big(Z_{n,k}- \Sbep[Z_{n,k} |\mathscr{H}_{n,k-1}]\big) Q_{k-1}^+.
\end{align*}
By the fact that $Z_{n,i}$s are bounded, Lemma \ref{lemma4.0.2} (f) and   H\"older's inequality, the  random variables considered above and in the sequel have finite moments of any order. So, the properties of the conditional expectation operator can be applied freely.  The sub-linear expectations of the last two sums above are non-positive, and   the sub-linear expectation of the second sum is also zero when   $\Sbep[Z_{n,k}|\mathscr{H}_{n,k}]\le 0$, $k=1,\ldots, k_n$. Taking the sub-linear expectation yields (\ref{eqlemRosenIeq.1}). By considering $\{-Z_{n,k}\}$, for $\max_{k\le k_n}(-S_{k_n}+S_k)$ we have a similar estimate.  Note $M_{k_n} \le 2 \max_{k\le k_{k_n}}|S_n-S_k|$. It follows that
\begin{align*}  \Sbep\left[M_{k_n}^2\right]\le & 8\Sbep\left[\sum_{k=1}^{k_n} \Sbep[Z_{n,k}^2|\mathscr{H}_{n,k-1}]\right]\\
& +16\Sbep\left[\sum_{k=1}^{k_n} \big\{(\Sbep[Z_{n,k}|\mathscr{H}_{n,k-1}])^+ + (\cSbep[Z_{n,k}|\mathscr{H}_{n,k-1}])^-\big\}M_{k_n}\right]\\
\le & 8\Sbep\left[\sum_{k=1}^{k_n} \Sbep[Z_{n,k}^2|\mathscr{H}_{n,k-1}]\right]+\frac{1}{2}\Sbep\left[M_{k_n}^2\right]\\
& +128\Sbep\left[\left(\sum_{k=1}^{k_n} \big\{(\Sbep[Z_{n,k}|\mathscr{H}_{n,k-1}])^+ + (\cSbep[Z_{n,k}|\mathscr{H}_{n,k-1}])^-\big\}\right)^2\right],
\end{align*}
where the last inequality is due to $ab\le \frac{a^2+b^2}{2}$.

For  (\ref{eqlemRosenIeq.4}), we apply the     elementary inequality
$$ |x+y|^p \le 2^{p}p^2|x|^p+|y|^p+px|y|^{p-1}\text{sgn} y+2^p p^2 x^2|y|^{p-2}, \;\;    p\ge 2, $$
and yields
$$ |Q_k|^p\le 2^pp^2|Z_{n,k}|^p+|Q_{k-1}|^p+pZ_{n,k}(Q_{k-1}^+)^{p-1}+2^p p^2Z_{n,k}^2(Q_{k-1}^+)^{p-2}. $$
It follows that
\begin{align*}
& \left(\max_{k\le k_n}(S_{k_n}-S_k)\right)^p\le    |Q_{k_n}|^p \\
\le & 2^pp^2\sum_{k=1}^{k_n} |Z_{n,k}|^p
+p\sum_{k=1}^{k_n}Z_{n,k}(Q_{k-1}^+)^{p-1}+2^p p^2\sum_{k=1}^{k_n}Z_{n,k}^2(Q_{k-1}^+)^{p-2}\\
\le & 2^pp^2\sum_{k=1}^{k_n} \Sbep[|Z_{n,k}|^p|\mathscr{H}_{n,k-1}]
+p\sum_{k=1}^{k_n}\big(\Sbep[Z_{n,k}|\mathscr{H}_{n,k-1}]\big)^+(Q_{k-1}^+)^{p-1}\\
& +2^p p^2\sum_{k=1}^{k_n}\Sbep[Z_{n,k}^2|\mathscr{H}_{n,k-1}](Q_{k-1}^+)^{p-2}
    + 2^pp^2\sum_{k=1}^{k_n} \left(|Z_{n,k}|^p-\Sbep[|Z_{n,k}|^p|\mathscr{H}_{n,k-1}]\right)\\
  &+p\sum_{k=1}^{k_n}\big(Z_{n,k}-\Sbep[Z_{n,k}|\mathscr{H}_{n,k-1}]\big) (Q_{k-1}^+)^{p-1}\\
  &+2^p p^2 \sum_{k=1}^{k_n}\big(Z_{n,k}^2-\Sbep[Z_{n,k}^2|\mathscr{H}_{n,k-1}]\big) (Q_{k-1}^+)^{p-2}.
\end{align*}
The sub-linear expectations of the last three sums are non-positive. Note $Q_k\le 2M_{k_n}$ and for $\big((\max_{k\le k_n} (-S_{k_n}+S_k)\big)^p$ we have a similar estimate. It follows that
\begin{align*}
\Sbep\left[M_{k_n}^p\right]\le & C_p\left\{\Sbep\left[\sum_{k=1}^{k_n} \Sbep[|Z_{n,k}|^p|\mathscr{H}_{n,k-1}]\right]
+\Sbep\left[\sum_{k=1}^{k_n} \Sbep[Z_{n,k}^2|\mathscr{H}_{n,k-1}] M_{k_n}^{p-2}\right]  \right.\\
& \left.+\Sbep\left[\sum_{k=1}^{k_n}\left\{\big(\Sbep[Z_{n,k}|\mathscr{H}_{n,k-1}]\big)^++\big(\cSbep[Z_{n,k}|\mathscr{H}_{n,k-1}]\big)^-\right\}M_{k_n}^{p-1}\right]
\right\}\\
\le & C_p\left\{\Sbep\left[\sum_{k=1}^{k_n} \Sbep[|Z_{n,k}|^p|\mathscr{H}_{n,k-1}]\right]
+\Sbep\left[\Big(\sum_{k=1}^{k_n} \Sbep[Z_{n,k}^2|\mathscr{H}_{n,k-1}]  \Big)^{p/2}\right]  \right.\\
& \left.+\Sbep\left[\Big(\sum_{k=1}^{k_n}\left\{\big(\Sbep[Z_{n,k}|\mathscr{H}_{n,k-1}]\big)^++\big(\cSbep[Z_{n,k}|\mathscr{H}_{n,k-1}]\big)^-\right\}\Big)^p
\right]
\right\}+\frac{1}{2} \Sbep[M_{k_n}^p],
\end{align*}
where the last inequality is due to $ab\le \frac{2}{p}|a|^{p/2}+(1-\frac{2}{p})|b|^{p/(p-2)}$ and $ab\le \frac{1}{p}|a|^p+(1-\frac{1}{p})|b|^{p/(p-1)}$. The proof is completed. $\Box$

\bigskip
The next one gives the exponential inequality of the martingale like sequences.
 \begin{theorem}\label{lem4.1}  Suppose that the operators $\Sbep_{n,k}$ satisfy (a) and (b), $\{Z_{n,k}; k=1,\ldots, k_n\}$ is an array of random variables such that $Z_{n,k}\in \mathscr{H}_{n,k}$ and $\Sbep[Z_{n,k}^2]<\infty$, $k=1,\ldots, k_n$.  Assume that $\Sbep[Z_{n,k}|\mathscr{H}_{n,k-1}]\le 0$ in $L_1$, $k=1,\ldots, k_n$.   Then for all $x,y,A>0$
 \begin{align}\label{eqlem4.1.1}
 \Capc\left(\max_{m\le k_n}\sum_{k=1}^{m} Z_{n,k}\ge x\right)\le &\Capc\left(\max_{k\le k_n} Z_{n,k}\ge y\; or \; \sum_{k=1}^{k_n} \Sbep[Z_{n,k}^2|\mathscr{H}_{n,k-1}]\ge A\right)\nonumber\\ &+ \exp\left\{-\frac{x^2}{2(xy+A)}\Big(1+\frac{2}{3}\ln \big(1+\frac{xy}{A}\big)\Big)\right\}.
 \end{align}
 \end{theorem}

 {\bf Proof.} Let $X_k=Z_{n,k}\wedge y$.  Then $Z_{n,k}-X_k=(Z_{n,k}-y)^+\ge 0$. Denote $\sigma_{n,k}^2=\Sbep[Z_{n,k}^2|\mathscr{H}_{n,k-1}]$, $\delta_k=\sum_{i=1}^k\sigma_{n,i}^2$, $k=1,\ldots, k_n$.
 Let $f(x)$ be a function with bounded derivative such that $I\{x\le  A\}\le f(x)\le I\{x\le A+ \epsilon\}$. Let $Y_k =X_kf(\delta_k)$, $T_k=\sum_{i=1}^k Y_k$. Then
 $\Sbep[Y_k|\mathscr{H}_{n,k-1}] \le  f(\delta_k)\Sbep[Z_{n,k}|\mathscr{H}_{n,k-1}]\le 0$    in $L_1$, $\Sbep[Y_k^2|\mathscr{H}_{n,k-1}] \le   f^2(\delta_k)\Sbep[Z_{n,k}^2|\mathscr{H}_{n,k-1}]= f^2(\delta_k)  \sigma_{n,k}^2$ in $L_1$. Denote $\delta^{\ast}_k =\sum_{i=1}^k f^2(\delta_k)  \sigma_{n,k}^2$.
 It follows that for any $x,y,A>0$,
 $$
 \Capc\left(\max_{m\le k_n}\sum_{k=1}^{m}Z_{n,k}\ge x\right)\le  \Capc\big(\max_{k\le k_n}Z_{n,k}\ge y\; or \;  \delta_{k_n}> A\big)+
 \Capc\left(\max_{k\le k_n} T_{k}\ge x\right).
 $$
  For any $t>0$, by noting $Y_k\le y$, $0\le f^2(\delta_k) \sigma_{n,k}^2\le \delta_k^{\ast}\le A+\epsilon$, and
$$e^{tY_k}=1+ tY_k+\frac{e^{tY_k}-1-t Y_k}{Y_k^2}Y_k^2\le 1 +tY_k+\frac{e^{ty}-1-t y}{y^2}Y_k^2, $$
we have
\begin{align*}
 &\exp\left\{-\frac{e^{ty}-1-t y}{y^2}f^2(\delta_k) \sigma_{n,k}^2
\right\}\Sbep\left[e^{t Y_k}\big|\mathscr{H}_{n,k-1}\right]\\
\le &\exp\left\{-\frac{e^{ty}-1-t y}{y^2}f^2(\delta_k) \sigma_{n,k}^2
\right\}\left\{ 1+\frac{e^{ty}-1-t y}{y^2} \Sbep[Y_k^2|\mathscr{H}_{n,k-1}]\right\} \\
\le &  1  \;\;\text{ in } L_1.
\end{align*}
Write
$$ U_0=1, \;\; U_k=\exp\Big\{-\frac{e^{ty}-1-t y}{y^2} \delta_k^{\ast} \Big\}e^{t T_k},\;\; k=1,\cdots, k_n. $$
Then
\begin{equation}\label{eqprooflem4.1.1} \Sbep\left[U_k|\mathscr{H}_{n,k-1}]\right]\le U_{k-1}\; \text{ in } L_1, \;\; k=1,\cdots, k_n.
\end{equation}
Next, we show that for any $\alpha>0$,
\begin{equation}\label{eqprooflem4.1.2}
\Capc\left(\max_{k\le k_n}U_k\ge \alpha\right)\le \frac{\Sbep[U_0]}{\alpha}.
\end{equation}
For given $\beta\in (0,\alpha)$, let $f(x)$ be a continuous function with bounded derivation such that $I\{x\le \alpha-\beta\}\le f(x) \le I\{x\le \alpha\}$. Define $f_0=1$, $f_k=f(U_1)\cdots f(U_k)$. Then $f_k\in \mathscr{H}_k$, $0\le f_k\le 1$ and
\begin{align*}
&f_0U_0+\sum_{k=1}^n f_{k-1}\big(U_k-U_{k-1}\big)=f_nU_n+\sum_{k=1}^n f_{k-1}\big(1-f(U_k)\big)U_k\\
\ge & f_nU_n+\sum_{k=1}^n f_{k-1}\big(1-f(U_k)\big)(\alpha-\beta)=(\alpha-\beta)(1-f_n)+f_nU_n\\
\ge & (\alpha-\beta)I\{\max_{k\le k_n}U_k\ge \alpha\}.
\end{align*}
By (\ref{eqprooflem4.1.1}),
\begin{align*} \Sbep\left[f_{k-1}\big(U_k-U_{k-1}\big)\right]= & \Sbep\left[\Sbep\left[f_{k-1}\big(U_k-U_{k-1}\big)\big|\mathscr{H}_{k-1}\right]\right]\\
=&\Sbep\left[f_{k-1}\big(\Sbep[U_k|\mathscr{H}_{k-1}]-U_{k-1}\big)\right]\le 0.
\end{align*}
It follows that
$$(\alpha-\beta)\Capc\left(\max_{k\le k_n}U_k\ge \alpha\right)\le \Sbep[f_0U_0]=\Sbep[U_0]. $$
(\ref{eqprooflem4.1.2}) is proved.
Now, note $\delta_k^{\ast}\le A+\epsilon$. We have for any $t>0$,
$$ \exp\left\{t\max_{k\le k_n}T_k\right\}\le \max_{k\le k_n}U_k \exp\Big\{\frac{e^{ty}-1-t y}{y^2} (A+\epsilon) \Big\}. $$
Hence by (\ref{eqprooflem4.1.2}),
\begin{align*}
\Capc\left(\max_{k\le k_n} T_k\ge x\right)  \le & \Capc\left(\max_{k\le k_n} U_k\ge \exp\Big\{tx-\frac{e^{ty}-1-t y}{y^2} (A+\epsilon) \Big\}\right)\\
 \le &  \exp\left\{-tx+\frac{e^{ty}-1-t y}{y^2} (A+\epsilon)\right\}.
 \end{align*}
Choosing $t=\frac{1}{y}\ln \big(1+\frac{xy}{A+\epsilon}\big)$ yields
$$
\Capc\left(\max_{k\le k_n} T_k\ge x\right)  \le  \exp\left\{\frac{x}{y}-\frac{x}{y}\Big(\frac{A+\epsilon}{xy}+1\Big)\ln\Big(1+\frac{xy}{A+\epsilon}\Big)\right\}.
$$
 Applying the elementary inequality
 $$ \ln (1+t)\ge \frac{t}{1+t}+\frac{t^2}{2(1+t)^2}\big(1+\frac{2}{3} \ln (1+t)\big)$$
 yields
 $$ \Big(\frac{A+\epsilon}{xy}+1\Big)\ln\Big(1+\frac{xy}{A+\epsilon}\Big)
 \ge 1+\frac{xy}{2(xy+A+\epsilon)}\Big(1+\frac{2}{3}\ln\big(1+\frac{xy}{A+\epsilon}\big)\Big). $$
 (\ref{eqlem4.1.1}) is proved by letting $\epsilon\to 0$. $\Box$

 \section{ L\'evy's characterization of a G-Brownian motion.}\label{sectLevy}
\setcounter{equation}{0}

In this section,    we   give a L\'evy   characterization of a G-Brownian motion as an application of Theorem \ref{thFCLTM}.    Let $\{\mathscr{H}_t; t\ge 0\}$   be  a non-decreasing family of subspaces of $\mathscr{H}$ such that (1) a constant $c\in \mathscr{H}_t$ and, (2) if $X_1,\ldots,X_d\in \mathscr{H}_t$, then $\varphi(X_1,\ldots,X_d)\in \mathscr{H}_t$ for any $\varphi\in C_{l,lip}$.
 We consider a system of operators in $\mathscr{L}(\mathscr{H})$,
 $$\Sbep_t: \mathscr{L}(\mathscr{H})\to \mathscr{L}(\mathscr{H}_t) $$
 and denote $\Sbep[X|\mathscr{H}_t]=\Sbep_t[X]$, $\cSbep[X|\mathscr{H}_t]=-\Sbep_t[-X]$.
Suppose that the operators $\Sbep_t$ satisfy  the following properties:  for all $X, Y \in \mathscr{L}({\mathscr H})$,
\begin{description}
  \item[\rm (i)]   $ \Sbep_t [ X+Y]=X+\Sbep_t[Y]$ in $L_1$ if $X\in \mathscr{H}_t$, and $ \Sbep_t [ XY]=X^+\Sbep_t[Y]+X^-\Sbep_t[-Y]$ in $L_1$
if $X\in \mathscr{H}_t$ and $XY\in \mathscr{L}({\mathscr H})$;
\item[\rm (ii)]   $\Sbep\left[\Sbep_t [ X]\right]=\Sbep[X]$.
 \end{description}

  \begin{example} Let $W_t$ be a G-Brownian motion in a sub-linear expectation space $(\Omega,\mathscr{H}, \Sbep)$, and
    $$\widetilde{\mathscr{H}}=\left\{X=\varphi(W_{t_1},\ldots, W_{t_d}): 0\le t_1\le \cdots\le t_d, \varphi\in C_{l,Lip}(\mathbb R_d), d\ge 1\right\}, $$
  $$\mathscr{H}_t=\left\{X=\varphi(W_{t_1},\ldots, W_{t_d}): 0\le t_1\le \cdots\le t_d\le t, \varphi\in C_{l,Lip}(\mathbb R_d), d\ge 1\right\}. $$
  For  $X= \varphi(W_{t_1},\cdots, W_{t_d})\in \widetilde{\mathscr{H}}$, assume $0\le t_1\le t_i\le t\le t_{i+1}\le \cdots\le t_d$, and define
  $$\Sbep_t[X]=\Sbep\left[\varphi(w_{t_1},\cdots,w_{t_i}, W_{t_{i+1}}-W_t+w_t,\cdots, W_{t_d}-W_t+w_t)\right]\Big|_{w_{t_1}=W_{t_1},\cdots,w_{t_i}=W_{t_i}, w_t=W_t}. $$
 Then, in the sub-linear expectation space  $(\Omega,\widetilde{\mathscr{H}}, \Sbep)$, the family  $\{\mathscr{H}_t, \Sbep_t\}_{t\ge 0}$ satisfies the properties (i)-(iii).
  \end{example}

 \begin{definition} A process $M_t$ is called a martingale, if $M_t\in \mathscr{L}(\mathscr{H})$, $M_t\in \mathscr{H}_t$ and
 $$ \Sbep[M_t|\mathscr{H}_s]=M_s, \;\; s\le t. $$
 \end{definition}
 Denote
  $$w_T(M,\delta)=\sup\limits_{|t-s|<\delta, t,s\in[0,T]}|M(t)-M(s)| $$
  and
 \begin{align*}
 W_T(M,\delta) = &\sup_{t_i}\Sbep\left[\max_{1\le i\le n} |M(t_i)-M(t_{i-1})|\wedge 1\right], \\
 & \text{ where the supermum }  \sup_{t_i} \text{ is taken over all } t_is \text{ with } \\
 & 0=t_0<t_1<\cdots<t_n=T,\;\; \delta/2<t_i-t_{i-1}<\delta, \; i=1,\cdots, n.
 \end{align*}

 The following theorem gives a L\'evy characterization of a G-Brownian motion.
  \begin{theorem} \label{th4.2} Let $M_t$ be a random process in $(\Omega,\mathscr{H},\mathscr{H}_t, \Sbep)$ with $M_0=0$,
 \begin{equation} \label{eqth4.2.1} \text{ for all }  p>0  \text{ and } t\ge 0, \;\; C_{\Capc}(|M_t|^p)<\infty \implies \Sbep[|M_t|^p]<\infty.
 \end{equation}    Suppose that $M_t$ satisfies
 \begin{description}
 \item[\rm (I)] both $M_t$ and $-M_t$ are martingales;
 \item[\rm (II)]  for a constant $\overline{\sigma}^2>0$, $M_t^2-\overline{\sigma}^2 t$ is a   martingale;
  \item[\rm (III)]
  for a constant $0<\underline{\sigma}^2\le \overline{\sigma}^2$,  $-(M_t^2-\underline{\sigma}^2 t)$ is a   martingale;
  \item[\rm (IV)]  for any $T>0$, $\lim_{\delta\to 0}W_T(M,\delta)=0$.
 \end{description}
 Then, $M_t$ satisfies Property (ii) as in Definition \ref{DefG-B}    with $M_1\sim N(0,[\underline{\sigma}^2,\overline{\sigma}^2])$.

 \end{theorem}

 \begin{remark} The assumption (I) implies that
 $\Sbep[M_t-M_s|\mathscr{H}_s]=\cSbep[M_t-M_s|\mathscr{H}_s]=0$ for all $t>s$. Also, under the assumptions (I),   the assumption (II)  is equivalent to  that
 $\Sbep[(M_t-M_s)^2|\mathscr{H}_s]=\overline{\sigma}^2(t-s)$ for all $t>s$,  (III)  is equivalent to that
 $\cSbep[(M_t-M_s)^2|\mathscr{H}_s]=\underline{\sigma}^2(t-s)$ for all $t>s$.

The assumption  of (IV)   means that $M_t$ is continuous.
Note $W_T(M,\delta)\le \epsilon+ \Capc\left(w_T(M,\delta)>\epsilon\right)$. It is satisfied if
\begin{description}
    \item[\rm (IV$^\prime$)] for any $T, \epsilon>0$, $  \lim_{\delta\to 0} \Capc\left( w_T(M,\delta)>\epsilon\right)=0.  $
 \end{description}
 The condition (IV$^{\prime}$) means that $M_t$ is continuous in capacity $\Capc$ uniformly in $t$ on each finite interval.
Also, $W_T(M,\delta)\le \sup_{t_i}\left( \sum_i \Sbep\left[|M(t_i)-M(t_{i-1})|^{2+\alpha}\right]\right)^{\frac{1}{2+\alpha}}$. (IV) is also satisfied if
 \begin{description}
 \item[\rm (IV)$^{\prime\prime}$] there is a constant  $\alpha>0$  such that  for any $t>s>0$, $\Sbep[|M_t-M_s|^{2+\alpha}]=o(t-s)$ as $t-s\to 0$.
 \end{description}
 \end{remark}

\begin{remark}
 The L\'evy characterization of a G-Brownian motion is first established under G-expectation   in a Wiener space by  Xu and Zhang (2009,2010) by using the stochastic calculus.  We will give an elementary proof by using the functional central limit theorem.
\end{remark}

  \begin{remark} If $\Sbep$ is countably sub-additive, then the condition (\ref{eqth4.2.1}) is automatically satisfied.
 The $G$-expectation space considered in Xu and Zhang (2009, 2010) is complete and so the sub-linear expectation is countably additive, and  (\ref{eqth4.2.1}) is satisfied.

In Xu and Zhang (2009, 2010),  the operators $\Sbep_t$ are also supposed to have the following assumptions:
\begin{description}
  \item[\rm (iii)] if $X\le Y$, then $\Sbep_t[X]\le   \Sbep_t[Y]$;
 \item[\rm (iv)]  $\Sbep_t[X]-\Sbep_t[Y]\le \Sbep_t[X-Y]$;
\item[\rm (v)]   $\Sbep_t\left[\left[\Sbep_s [ X]\right]\right]=\Sbep_{t\wedge s} [ X]$.
 \end{description}
 As in Lemma \ref{lemma4.0.2}, (iii), (iv), and (v) holds in $L_1$ if the operators satisfy (i) and (ii).
 \end{remark}

For proving Theorem \ref{th4.2} we need   a more lemma.

 \begin{lemma}\label{lem4.2}  Suppose that the operators $\Sbep_t$ satisfy (i)-(ii), $M_t$ is a martingale in $(\Omega, \mathscr{H}, \mathscr{H}_t, \Sbep)$ such that (IV) in Theorem \ref{th4.2} is satisfied and $\Sbep[(M_t-M_s)^2|\mathscr{H}_s]\le   (t-s)\sigma^2$ for all $t>s\ge 0$, where $\sigma$ is a positive constant.   Then,
  \begin{equation} \label{eqlem4.2.1}
  \Capc\left( M_t-M_s\ge x\right)\le \exp\left\{-\frac{x^2}{2   (t-s)\sigma^2}\right\}, \text{ for all } t>s\ge 0, x\ge 0.
  \end{equation}
  In particular, for any $p>0$, $C_{\Capc}\left(\big[(M_t-M_s)^+\big]^p\right)
   \le   c_p   (t-s)^{p/2}\sigma^p.$
 \end{lemma}

 {\bf Proof.} Let $s=t_0<t_1<\ldots<t_k=t$ be a partition of $[s,t]$ with $\delta/2<t_i-t_{i-1}<\delta$. Note
  $\Sbep[M_{t_i}-M_{t_{i-1}}|\mathscr{H}_{t_{i-1}}]=0$ and $\Sbep[(M_{t_i}-M_{t_{i-1}})^2|\mathscr{H}_{t_{i-1}}]\le     (t_i-t_{i-1})\sigma^2$. So, $\sum_{i=1}^k \Sbep[(M_{t_i}-M_{t_{i-1}})^2|\mathscr{H}_{t_{i-1}}]\le    (t-s)\sigma^2$.  By Theorem \ref{lem4.1}, for $0<y<1$ and $x>0$,
  \begin{align*}
  &\Capc\left( M_t-M_s\ge x\right)\\
  \le &  \Capc\left(\max_i (M_{t_i}-M_{t_{i-1}})\ge y\right) \\
  & \;\; +\exp\left\{-\frac{x^2}{2(xy+    (t-s)\sigma^2)}\Big(1+\frac{2}{3}\ln \big(1+\frac{xy}{   (t-s)\sigma^2}\big)\Big)\right\}
  \\
  \le &  \frac{W_T(M,\delta)}{y} +\exp\left\{-\frac{x^2}{2(xy+   (t-s)\sigma^2)}\Big(1+\frac{2}{3}\ln \big(1+\frac{xy}{  (t-s)\sigma^2}\big)\Big)\right\}.
  \end{align*}
 By letting $\delta\to 0$ and then $y\to 0$, we conclude  (\ref{eqlem4.2.1}). Finally, for  $p>0$,
   \begin{align*} & C_{\Capc}\left(\big[(M_t-M_s)^+\big]^p\right)\le   \int_0^{\infty}\Capc\left( M_t-M_s \ge x^{1/p}\right)dx \\
   \le  &   (t-s)^{p/2}\sigma^p\int_0^{\infty}\exp\left\{-\frac{x^{2/p}}{2}\right\}dx
   \le   c_p  (t-s)^{p/2}\sigma^p. \qquad \Box
   \end{align*}
 \bigskip

 {\bf Proof of Theorem \ref{th4.2}}.  Suppose that (I)-(IV) are satisfied.   Note that both $M_t$ and $-M_t$ are martingales, and  $\Sbep[(M_t-M_s)^2|\mathscr{H}_{t_{i-1}}]=  (t-s)\overline{\sigma}^2$.    By Lemma \ref{lem4.2},
   \begin{align*} C_{\Capc}\left(|M_t-M_s|^p\right)
   \le   c_p  (t-s)^{p/2}\overline{\sigma}^p.
   \end{align*}
  By the assumption (\ref{eqth4.2.1}), $\Sbep[|M_t-M_s|^p]<\infty$ for any $p>0$ and $t,s$. Let $W_t$ be a G-Brownian motion in a sub-linear expectation $(\widetilde{\Omega}, \widetilde{\mathscr{H}}, \widetilde{\mathbb E})$ with $W_1\sim N(0,[\underline{\sigma}^2,\overline{\sigma}^2])$.  It is sufficient to show that for any  $0<t_1<\ldots<t_d$  and   $\varphi\in C_{b,Lip}(\mathbb R_d)$
 \begin{equation}\label{eqpproofth4.2.1} \Sbep\left[\varphi(M_{t_1},\cdots, M_{t_d})\right]=\widetilde{\mathbb{E}}\left[\varphi(W_{t_1},\cdots, W_{t_d})\right].
 \end{equation}
   Actually, by noting $\Sbep[|M_t|^p]<\infty$ for any $p>0$, we can extend $\varphi$ from $C_{b,Lip}(\mathbb R_d)$ to $C_{l,Lip}(\mathbb R_d)$ by an elementary argument.

 Now, without loss of generality, we assume $0<t_1<\ldots<t_d\le 1$.  Note $\Sbep\left[(|M_t-M_s|^3-c^3)^+\right]\le \Sbep\left[|M_t-M_s|^4\right]/c \to 0$ as $c\to \infty$. Then  $\Sbep[|M_t-M_s|^3]\le C_{\Capc}\left(|M_t-M_s|^3\right)=o(t-s)$ as $t-s\to 0$. Let
 $$ k_n=2^n, \;\; Z_{n,k}=M_{k/2^n}-M_{(k-1)/2^n}, \;\; \mathscr{H}_{n,k}=\mathscr{H}_{k/2^n}, \;\; k=1,\ldots, k_n, $$
 and  $\tau_n(t)=[t2^n]$. Then $\Sbep[Z_{n,k}|\mathscr{H}_{n,k-1}]=\cSbep[Z_{n,k}|\mathscr{H}_{n,k-1}]=0$,
 $$\Sbep[Z_{n,k}^2|\mathscr{H}_{n,k-1}]=\frac{\overline{\sigma}^2}{2^n}, \;\; \cSbep[Z_{n,k}^2|\mathscr{H}_{n,k-1}]=\frac{\underline{\sigma}^2}{2^n}. $$
 Hence it is easily seen that the sequence $\{Z_{n,k}, \mathscr{H}_{n,k}\}$ satisfy the conditions (\ref{eqCLTCondM.3}), (\ref{eqCLTCondM.4}) and (\ref{eqFCLTCondM.2}) with $\rho(t)=t\overline{\sigma}^2$, $r=\underline{\sigma}^2/\overline{\sigma}^2$. Further,
 $$ \sum_{k=1}^{k_n}\Sbep[|Z_{n,k}|^3]=\sum_{k=1}^{2^n}o\big(\frac{1}{2^n}\big)\to 0. $$
 So, the Lindeberg condition (\ref{eqLindebergM}) is satisfied. Let $W_n(\cdot)$ be defined as in (\ref{eqthFCLTM.1}). By Theorem \ref{thFCLTM},
 $  (W_n(t_1),\cdots,W_n(t_d))\overset{d}\to (W_{t_1},\cdots, W_{t_d}). $
 On the other hand,
 $$|W_n(t)-M_t|= \Big|M_t-M_{[2^nt]/2^n}\Big|\overset{\Capc}\to 0. $$
 So, (\ref{eqpproofth4.2.1}) holds for all $\varphi\in C_{b,Lip}(\mathbb R_d)$. The proof is now completed.   $\Box$

 \bigskip
\section{Proofs of the central limit theorems for martingales.}\label{sectProof}
  \setcounter{equation}{0}

\subsection{Proof of the central limit theorem}

We give the proof of Theorem \ref{thCLTM}.   By (\ref{eqLindebergM}),   there exists a sequence of positive numbers  $1/2>\epsilon_n\searrow 0$ such that
$$\epsilon_n^{-2}  \sum_{k=1}^{k_n}\Sbep\left[ \left( Z_{n,k}^2-\epsilon_n^2 \right)^+|\mathscr{H}_{n,k-1}\right]  \overset{\Capc}\to 0. $$
Let $Z_{n,k}^{\ast}=(-2\epsilon_n)\vee Z_{n,k}\wedge (2\epsilon_n)$. Then
$$    \sum_{k=1}^{k_n}\Sbep\left[ \left( Z_{n,k}-Z_{n,k}^{\ast} \right)^2|\mathscr{H}_{n,k-1}\right]  \le  \sum_{k=1}^{k_n}\Sbep\left[ \left( Z_{n,k}^2-\epsilon_n^2 \right)^+|\mathscr{H}_{n,k-1}\right] \overset{\Capc}\to 0$$
and
$$ \sum_{k=1}^{k_n}\Sbep\left[  | Z_{n,k}-Z_{n,k}^{\ast}  | \big|\mathscr{H}_{n,k-1}\right] \le \epsilon_n^{-1}  \sum_{k=1}^{k_n}\Sbep\left[ \left( Z_{n,k}^2-\epsilon_n^2 \right)^+|\mathscr{H}_{n,k-1}\right] \overset{\Capc} \to 0.$$
Hence, $\{Z_{n,k}^{\ast}; k=1,\ldots, k_n\}$ satisfy the conditions (\ref{eqCLTCondM.2})-(\ref{eqCLTCondM.4}). Further, let
$h_k=\epsilon_n^{-2}\sum_{i=1}^k \Sbep\big[ \big( Z_{n,k}^2-\epsilon_n^2 \big)^+|\mathscr{H}_{n,k-1}\big]$ and $f$ be a bounded Lipschitz function  such that $I\{x\le \epsilon\}\le f(x)\le I\{x\le 2\epsilon\}$. Then,
\begin{align*}
&\Capc\left( Z_{n,k}\ne Z_{n,k}^{\ast} \text{ for some } k\right)\\
 = &\Capc\left( \max_{k\le k_n}|Z_{n,k}|\ge 2\epsilon_n\right)
\le   \Capc\left(\sum_{k=1}^{k_n}\left[1\wedge (Z_{n,k}^2-\epsilon_n^2)^+\right]\ge \epsilon_n^2\right) \\
\le &  \Capc\left(\sum_{k=1}^{k_n} (\left[1\wedge (Z_{n,k}^2-\epsilon_n^2)^+\right]\ge \epsilon_n^2, h_{k_n}\le \epsilon\right)+\Capc(h_{k_n}\ge \epsilon)\\
= &  \Capc\left(\sum_{k=1}^{k_n} (\left[1\wedge (Z_{n,k}^2-\epsilon_n^2)^+\right]f(h_k)\ge \epsilon_n^2, h_{k_n}\le \epsilon\right)+\Capc(h_{k_n}\ge \epsilon)\\
\le & \Sbep\left[\epsilon_n^{-2}\sum_{k=1}^{k_n} \left[1\wedge (Z_{n,k}^2-\epsilon_n^2)^+\right]f(h_k)\right]+\Capc(h_{k_n}\ge \epsilon)\\
\le & \Sbep\left[\epsilon_n^{-2}\sum_{k=1}^{k_n}f(h_k)\Sbep\left[ \left[1\wedge (Z_{n,k}^2-\epsilon_n^2)^+\right]|\mathscr{H}_{n,k-1}\right] \right]+\Capc(h_{k_n}\ge \epsilon)\\
\le &2\epsilon+\Capc(h_{k_n}\ge \epsilon)\to 0 \text{ as } n\to \infty \text{ and then } \epsilon\to 0.
\end{align*}
It follows that for any bounded function $\varphi$,
$$ \Sbep\Big[\big|\varphi(\sum_{k=1}^{k_n}Z_{n,k})-\varphi(\sum_{k=1}^{k_n}Z_{n,k}^{\ast} )\big|\Big]
\le 2\sup_x|\varphi(x)|\Capc \left( Z_{n,k}\ne Z_{n,k}^{\ast} \text{ for some } k\right)\to 0. $$
So, without loss of generality we can assume that there is a positive sequence $1\ge \epsilon_n\searrow 0$ such that $|Z_{n,k}|\le \epsilon_n$, $k=1,\ldots, k_n$.

Denote  $S_0=0$, $\delta_0=0$, $S_k=\sum_{i=1}^k Z_{n,i}$, $a_{n,k}^2=\Sbep[Z_{n,k}^2|\mathscr{H}_{n,k-1}]$, $\delta_k=\sum_{i=1}^ka_{n,i}^2$, $k=1,\ldots, k_n$.
 Let $f(x)$ be a function with bounded derivative such that $I\{x\le \rho+\epsilon/2\}\le f(x)\le I\{x\le \rho+ \epsilon\}$. Let $Z_{n,k}^{\ast}=Z_{n,k}f(\delta_k)$. Then $\{Z_{n,k}^{\ast}; k=1,\ldots, k_n\}$ satisfy the conditions (\ref{eqCLTCondM.2})-(\ref{eqCLTCondM.4}), and
 \begin{equation}\label{eqProofCLTC0.1}  \sum_{k=1}^{k_n} \Sbep[(Z_{n,k}^{\ast})^2|\mathscr{H}_{n,k-1}]=\delta_{k_n}^{\ast},
 \end{equation}
where  $\delta_{k_n}^{\ast}=
  \sum_{k=1}^{k_n} f(\delta_k)\Sbep[(Z_{n,k}^2|\mathscr{H}_{n,k-1}]\le  \rho+\epsilon. $
  The above  equalities  hold in $L_1$ by the Property (a) of the operators  $\Sbep_{n,k}$ and then hold  in any $L_q$ by Lemma \ref{lemma4.0.1} (2) since $\delta_{k_n}^{\ast}$ is bounded in $L_q$ by Lemma \ref{lemma4.0.2} (f).
Further,
$$\left\{  Z_{n,k} \ne   Z_{n,k}^{\ast}  \text{ for some } k\right\}\subset \left\{ \sum_{k=1}^{k_n} a_{n,k}^2>\rho+\epsilon/2\right\}. $$
 So, without loss of generality we can further assume that
$
  \delta_{k_n}=\sum_{k=1}^{k_n}\Sbep[Z_{n,k}^2|\mathscr{H}_{n,k-1}]\le \rho+\epsilon$ in $L_1$.
   Similarly, we can assume
$
 \chi_{k_n}=:\sum_{k=1}^{k_n}\left\{ |\Sbep[Z_{n,k} |\mathscr{H}_{n,k-1}]|+|\cSbep[Z_{n,k} |\mathscr{H}_{n,k-1}]|\right\}<\epsilon<1
$ in $L_1$.
The property (f) in Lemma \ref{lemma4.0.2} implies  that  all random variables considered above and in the sequel are bounded in $L_p$ for all $p>0$.

Now, by Theorem \ref{lemRosenIeq},
 \begin{equation}\label{eqvariance}
 \Sbep\left[\max_{k\le k_n} \big(\sum_{i=1}^{k}Z_{n,i} \big)^2\right]
 \le 256 \Sbep\left[ \delta_{k_n}\right] +
  256\Sbep\left[ \chi_{k_n}^2 \right].
 \end{equation}
  If $\rho=0$, then $\delta_{k_n}\overset{\Capc}\to 0$. Note $\chi_{k_n}\overset{\Capc}\to 0$. So,  $\Sbep\left[  \big(\sum_{i=1}^{k_n}Z_{n,i} \big)^2\right]\to 0$, and then    the result is obvious.
 In the sequel, we suppose $\rho\ne 0$.  Let  $\varphi$ be a  bounded  continuous function with bounded derivation. Without loss of generality, we assume $|\varphi(x)|\le 1$.
 We want to show that
 \begin{equation}\label{eqconvergence}\Sbep[\varphi(S_{k_n})]\to \widetilde{\mathbb E}[\varphi(\sqrt{\rho}\xi)].
 \end{equation}
 In the classical probability space, the above convergence is usually shown by verifying  the convergence of the related characteristic functions
 (cf. Hall and Heyde (1980), p. 60-63; Pollard (1984), p. 171-174). As shown by  Hu and Li (2014),   the characteristic function cannot determine the distribution of
random variables in the sub-linear expectation space.   Peng (2007a, 2008b) developed a   method to show the above convergence for independent random variables. Here we promote Peng's argument such that it is also valid for martingale differences which give also a new normal approximation method for classical martingale differences instead of the characteristic function.

 Now, for a small but fixed $h > 0$, let $V (t, x)$ be the unique viscosity solution of the following equation,
\begin{equation}\label{eqPDE} \partial_t V + G( \partial_{xx}^2 V)=0,\;\;  (t, x) \in [0,\rho+ h] \times \mathbb R, \; V|_{t=\rho+h} = \varphi(x),
\end{equation}
where $G(\alpha)=\frac{1}{2}\big( \alpha^+-r\alpha^-\big)$. Then by
the interior regularity of $V$,
\begin{equation}\label{eqPDE2}\|V\|_{C^{1+\alpha/2,2+\alpha}([0,\rho+h/2]\times R)} < \infty, \text{ for some } \alpha\in (0, 1).
\end{equation}
According to the definition of $G$-normal distribution, we have $V(t,x)=\widetilde{\mathbb E}\big[\varphi(x+\sqrt{\rho+h-t}\xi)\big]$ where $\xi\sim N(0,[r,1])$ under $\widetilde{\mathbb E}$. In particular,
$$ V(h,0)=\widetilde{\mathbb E}\big[\varphi( \sqrt{\rho}\xi)\big], \;\; V(\rho+h, x)=\varphi(x). $$
It is obvious that, if  $\varphi(\cdot)$ is a global Lipschitz function, i.e.,
$|\varphi(x)-\varphi(y)|\le C|x-y|$, then $|V(t,x)-V(t,y)|\le C|x-y|$ and
$$|V(t,x)-V(s,x)|\le C\widetilde{\mathbb E}[|\xi|] \left|\sqrt{\rho+h-t}-\sqrt{\rho+h-s}\right| \le C\widetilde{\mathbb E}[|\xi|] |t-s|^{1/2}.$$
So, $|\partial_x V(t,x)|\le C$, $|\partial_t V(t,x)|\le  C\widetilde{\mathbb E}[|\xi|]/\sqrt{\rho+h-t}$,  $|V(\rho+h,x)-V(\rho,x)|\le C\widetilde{\mathbb E}[|\xi|] \sqrt{h}$ and $|V(h,0)-V(0,0)|\le C\widetilde{\mathbb E}[|\xi|] \sqrt{h}$.    Following the proof of Lemma 5.4 of Peng (2008b), it is sufficient to show that
\begin{equation}\label{eqconV}
 \lim_{n\to \infty} \Sbep[ V(\rho, S_{k_n})]=V(0,0).
 \end{equation}
 As we have shown, we can assume that $\delta_{k_n} \le \rho+h/4 =:h_0<2\rho$ in $L_1$.  It is obvious that $|V(t,x)|\le 1$, and
$$ \Sbep\left[\Big|V(\rho,S_{k_n})-V(\delta_{k_n}\wedge h_0,S_{k_n})\Big|\right]\le C\Sbep\left[ |\delta_{k_n}\wedge h_0-\rho|^{1/2}\right]\to 0. $$
Hence, it is sufficient to show that
\begin{equation}\label{eqtheorem2proof.1}
\lim_{n\to \infty} \Sbep\left[V(\delta_{k_n}\wedge h_0, S_{k_n})\right]=V(0,0).
\end{equation}
Let $\widetilde{\delta}_i =\delta_i \wedge h_0$.  Then $\widetilde{\delta}_{i+1}-\widetilde{\delta}_i\le a_{n,i+1}^2$, $|\widetilde{\delta}_i|\le h_0=\rho+h/4$. It follows that
$$ \Big|\partial_x V(\widetilde{\delta}_i,  S_i)\Big|\le C, \;\; \Big|\partial_t  V(\widetilde{\delta}_i,  S_i)\Big|\le C/\sqrt{h}\le C. $$
Also, by the fact that $\partial_{xx} V$ is uniformly $\alpha$-H\"older continuous in $x$ and $\alpha/2$-H\"older   continuous in t on $[0, \rho+h/2] \times R$, it follows that
$$ \Big|\partial_{xx}^2 V(\widetilde{\delta}_i,  S_i)\Big|\le  \Big|\partial_{xx}^2 V(0,0)\Big| +C|\widetilde{\delta}_i|^{\alpha/2}+C|S_i|^{\alpha}\le C+C|S_i|^{\alpha}. $$
Now,
applying the Taylor's expansion yields
\begin{align*}
 & V(\widetilde{\delta}_{k_n}, S_{k_n})-V(0,0)\\
 =&\sum_{i=0}^{k_n-1} \left\{[ V( \widetilde{\delta}_{i+1},  S_{i+1})-V( \widetilde{\delta}_i,  S_{i+1})]
  +[ V( \widetilde{\delta}_i,  S_{i+1})-V( \widetilde{\delta}_i,  S_i)]\right\}
 =: \sum_{i=0}^{k_n-1} \left\{I_{n}^i+J_{n}^i\right\},
\end{align*}
with
\begin{align*}
 J_{n}^i =&  \partial_t V( \widetilde{\delta}_i,  S_i)\big(\widetilde{\delta}_{i+1}-\widetilde{\delta}_i\big)    +\frac{1}{2} \partial_{xx}^2 V(\widetilde{\delta}_i,  S_i)Z_{n,i+1}^2    +
\partial_x V( \widetilde{\delta}_i,  S_i) Z_{n,i+1}       \\
=&      \Big\{a_{n,i+1}^2 \partial_t V( \widetilde{\delta}_i,  S_i)    +\frac{1}{2} \partial_{xx}^2 V(\widetilde{\delta}_i,  S_i)  Z_{n,i+1}^2    -\frac{1}{2}(\partial_{xx}^2 V(\widetilde{\delta}_i,  S_i))^- (ra_{n,i+1}^2-\cSbep[Z_{n,i+1}^2|\mathscr{H}_{n,i}])\Big\}  \\
&  + \Big\{ \partial_x V( \widetilde{\delta}_i,  S_i) Z_{n,i+1}\Big\}
   +\Big\{\frac{1}{2}(\partial_{xx}^2 V(\widetilde{\delta}_i,  S_i))^- (ra_{n,i+1}^2-\cSbep[Z_{n,i+1}^2|\mathscr{H}_{n,i}])\Big\}\\
 & +\Big\{\partial_t V( \widetilde{\delta}_i,  S_i)\big(\widetilde{\delta}_{i+1}-\widetilde{\delta}_i-a_{n,i+1}^2\big)\Big\}
   \\
 =:& J_{n,1}^i+J_{n,2}^i+J_{n,3}^i+J_{n,4}^i
\end{align*}
 and
\begin{align*}
I_{n}^i=&\big(\widetilde{\delta}_{i+1}-\widetilde{\delta}_i\big)\left[  \big( \partial_t V( \widetilde{\delta}_i+\gamma \big(\widetilde{\delta}_{i+1}-\widetilde{\delta}_i\big),  S_{i+1})-\partial_t V( \widetilde{\delta}_i ,  S_{i+1})\big) \right.\\
&\left.\qquad \quad +\big( \partial_t V( \widetilde{\delta}_i,  S_{i+1})-\partial_t V( \widetilde{\delta}_i ,  S_i)\big)   \right] \\
&+ \frac{1}{2}\left[ \partial_{xx}^2 V(\widetilde{\delta}_i,  S_i+\beta Z_{n,i+1})-\partial_{xx}^2 V(\widetilde{\delta}_i,  S_i)\right]Z_{n,i+1}^2,
\end{align*}
where $\gamma$ and $\beta$ are between $0$ and $1$.  Thus
\begin{align}\label{eqtheorem2proof.2}
 &\Big|\Sbep [V(\widetilde{\delta}_{k_n}, S_{k_n})]-V(0,0)-\Sbep\big[\sum_{i=0}^{k_n-1}(J_{n,1}^i+J_{n,2}^i)\big]\Big| \nonumber\\
  \le &
 \Sbep\Big[\big|V(\widetilde{\delta}_{k_n}, S_{k_n}) -V(0,0)-\sum_{i=0}^{k_n-1}(J_{n,1}^i+J_{n,2}^i)\big|\Big]\\
 \le &  \Sbep\Big[  \sum_{i=0}^{k_n-1}(|I_{n}^i|+|J_{n,3}^i|+|J_{n,4}^i|)\Big].\nonumber
\end{align}
For $J_{n,1}^i$,   it follows that
$$ \Sbep\left[J_{n,1}^i\big|\mathscr{H}_{n,i}\right]=
\big[\partial_t V(\widetilde{\delta}_i,  S_i)+G\big(\partial^2_{xx}   V( \widetilde{\delta}_i,  S_i)\big)\big]a_{n,i+1}^2  =0\;\text{ in } L_1. $$
It follows that \begin{equation}\label{eqtheorem2proof.3}\Sbep\left[\sum_{i=0}^{k_n-1}J_{n,1}^i\right]=\Sbep\left[\sum_{i=0}^{k_n-2}J_{n,1}^i+\Sbep\left[J_{n,1}^{k_n-1}\big|\mathscr{H}_{n,k_n-1}\right]\right]=
\Sbep\left[\sum_{i=0}^{k_n-2}J_{n,1}^i\right]=\ldots=0.
\end{equation}
For $J_{n,2}^i$, we denote $ \widetilde{J}_{n,2}^i=\big|\partial_x V( \widetilde{\delta}_i,  S_i)\big|\left(\big|\Sbep[ Z_{n,i+1}|\mathscr{H}_{n,i}]\big| + \big|\cSbep[ Z_{n,i+1}|\mathscr{H}_{n,i}]\big|\right)$. Then
\begin{align*}
&\Sbep[J_{n,2}^i-\widetilde{J}_{n,2}^i|\mathscr{H}_{n,i}]=\Sbep[J_{n,2}^i|\mathscr{H}_{n,i}]-\widetilde{J}_{n,2}^i\\
\le &  (\partial_x V( \widetilde{\delta}_i,  S_i))^+\Sbep[ Z_{n,i+1}|\mathscr{H}_{n,i}] - (\partial_x V( \widetilde{\delta}_i,  S_i))^-\cSbep[ Z_{n,i+1}|\mathscr{H}_{n,i}] -\widetilde{J}_{n,2}^i
\le   0 \text{ in } L_1.
\end{align*}
Similarly
$\Sbep[-J_{n,2}^i-\widetilde{J}_{n,2}^i|\mathscr{H}_{n,i}]\le 0$ in  $L_1$.
It follows that
 \begin{align}\label{eqtheorem2proof.4}  \Sbep\Big[\sum_{i=0}^{k_n-1}(\pm J_{n,2}^i-\widetilde{J}_{n,2}^i)\Big]
=&\Sbep\left[\sum_{i=0}^{k_n-2}(\pm J_{n,2}^i-\widetilde{J}_{n,2}^i)+\Sbep\left[\pm J_{n,1}^{k_n-1}-\widetilde{J}_{n,2}^{k_n-1}\big|\mathscr{H}_{n,k_n-1}\right]\right]\nonumber\\
\le & \Sbep\Big[\sum_{i=0}^{k_n-2}(\pm J_{n,2}^i-\widetilde{J}_{n,2}^i)\Big]\le \ldots\le 0.
\end{align}
Hence
\begin{equation}\label{eqtheorem2proof.5}
\Sbep\Big[\pm \sum_{i=0}^{k_n-1}J_{n,2}^i\Big] \le \Sbep\Big[\sum_{i=0}^{k_n-1}(\pm J_{n,2}^i-\widetilde{J}_{n,2}^i)\Big]
+\Sbep\Big[\sum_{i=0}^{k_n-1} \widetilde{J}_{n,2}^i\Big] \le \Sbep\Big[\sum_{i=0}^{k_n-1} \widetilde{J}_{n,2}^i\Big].
\end{equation}
Note $\big|\partial_x V( \widetilde{\delta}_i,  S_i)\big|\le C$,   $\chi_{k_n}\overset{\Capc}\to 0 $ and $\chi_{k_n}\le 1$  in any $L_p$. Combining (\ref{eqtheorem2proof.3}) and (\ref{eqtheorem2proof.5}) yields that
$$\Big|\Sbep\big[\sum_{i=0}^{k_n-1}(J_{n,1}^i+J_{n,2}^i)\big]\Big|\le \Sbep\Big[\sum_{i=0}^{k_n-1} \widetilde{J}_{n,2}^i\Big]
\le C\Sbep[\chi_{k_n}]\to 0. $$

For  $J_{n,3}^i$, it is easily seen that
\begin{equation}\label{eqProofCLTC0.4}
 \sum_{i=0}^{k_n-1} |J_{n,3}^i| \le C (1+\max_{i\le k_n}|S_i|^{\alpha}) \sum_{i=1}^{k_n}  |ra_{n,i}^2-\cSbep[Z_{n,i}^2|\mathscr{H}_{n,i-1}]|.
 \end{equation}
 Write $\beta_{k_n}= \sum_{i=1}^{k_n} |ra_{n,i}^2-\cSbep[Z_{n,i}^2|\mathscr{H}_{n,i-1}]|$.  Note that
  $$\beta_{k_n} \overset{\Capc}\to 0 \; \text{ and }\;  \beta_{k_n}  \le  2\delta_{k_n} \le 2h_0  \text{ in any } L_p$$
  and
   $\Sbep[\max_{i\le k_n}|S_i|^2]\le 256\{\Sbep[\delta_{k_n}]+\Sbep[\chi_{k_n}^2]\}\le 256(h_0+1) $ by (\ref{eqvariance}).
So
$$\Sbep\Big[  \sum_{i=0}^{k_n-1}  |J_{n,3}^i|  \Big] \le C\big(\Sbep[( 1+\max_{i\le k_n}|S_i|^{\alpha})^2]\big)^{1/2}\big(\Sbep[\beta_{k_n}^2]\big)^{1/2} \to 0.
$$

For  $J_{n,4}^i$, note that $\big|\widetilde{\delta}_{i+1}-\widetilde{\delta}_i-a_{n,i+1}^2\big|\le a_{n,i+1}^2$, and $\widetilde{\delta}_{i+1}-\widetilde{\delta}_i-a_{n,i+1}^2=\delta_{i+1}- \delta_i-a_{n,i+1}^2=0$
when $\delta_{k_n}\le h_0$. It follows that
$$\Sbep\left[\sum_{i=0}^{k_n-1}  |J_{n,4}^i|\right]\le C \Sbep\left[\delta_{k_n}I\{\delta_{k_n}>h_0\}\right]\le C  \left(\Sbep\left[\delta_{k_n}^2\right]\right)^{1/2}
\left(\Capc(\delta_{k_n}>h_0)\right)^{1/2}=0. $$

For $I_n^i$, note both $\partial_t V$ and $\partial_{xx} V$ are uniformly $\alpha$-H\"older continuous in $x$ and $\alpha/2$-H\"older   continuous in t on $[0, \rho+h/2] \times R$. Without loss of generality, we assume $\alpha<\tau$. Also, $\widetilde{\delta}_{i+1}-\widetilde{\delta}_i\le a_{n,i+1}$.  We then have
\begin{align*}
|I_n^i|\le  &  C\big|a_{n,i+1} \big|^{2+\alpha } +C a_{n,i+1}^2 |Z_{n,i+1}|^{\alpha} +|Z_{n,i+1}|^{2+\alpha}\\
\le &C\epsilon_n^{\alpha}a_{n,i+1}^2+C \epsilon_n^{\alpha} Z_{n,i+1}^2
=   C\epsilon_n^{\alpha}a_{n,i+1}^2+C\epsilon_n^{\alpha}\left(Z_{n,i+1}^2- a_{n,i+1}^2\right)
\end{align*}
  in any $L_q$ by Lemma \ref{lemma4.0.1}.
And so,
 \begin{equation} \label{eqProofCLTC0.5} \sum_{i=0}^{k_n-1}|I_{n}^i|\le 2 C\epsilon_n^{\alpha}+C\epsilon_n^{\alpha}\sum_{i=1}^{k_n}\left(Z_{n,i}^2- a_{n,i}^2\right) \;\;\text{ in } L_1,
 \end{equation}
by noting $ \sum_{i=1}^{k_n}a_{n,i}^2\le 2\rho$ in $L_1$, where the sub-linear expectation  under $\Sbep$ of the last   term is zero. It follows  that
$$
 \Sbep\Big[  \sum_{i=0}^{k_n-1}|I_{n}^i|\Big]\le 2 C\epsilon_n^{\alpha} \to 0.
$$
(\ref{eqtheorem2proof.1}) is proved. Hence, (\ref{eqconvergence}) holds for any bounded function $\varphi$ with bounded derivative.

If $\varphi$ is a bounded and  uniformly continuous function,
 we   define a function $\varphi_{\delta}$ as a convolution of $\varphi$ and the density of a normal distribution
$N(0,\delta)$, i.e.,
$$ \varphi_{\delta}= \varphi\ast \psi_{\delta}, \;\; \text{with} \;
\psi_{\delta}(x)=\frac{1}{\sqrt{2\pi\delta}}\exp\left\{-\frac{x^2}{2\delta}\right\}, $$
where $\varphi\ast \psi_{\delta}$ denotes the convolution of $\varphi $ and $\psi_{\delta}$.  Then $|\varphi_{\delta}^{\prime}(x)|\le \sup_x|\varphi(x)|\delta^{-1/2}$ and $\sup_x|\varphi_{\delta}(x)-\varphi(x)|\to 0$ as $\delta\to 0$.
 Hence, (\ref{eqconvergence}) holds for  any  bounded and uniformly continuous function $\varphi$.

 Now,   for a bounded continuous function  $\varphi$ and a give a number $N>1$, we
define $\varphi_1(x)=\varphi\big((-N)\vee (x\wedge N)\big)$. Then, $\varphi_1$ is a bounded and uniformly continuous function, and $|\varphi(x)-\varphi_1(x)|\le C I\{|x|>N\}$. And so,
\begin{align*}
& \sup_n\Sbep\Big[\big|\varphi(\sum_{k=1}^{k_n}Z_{n,k} )-\varphi_1(\sum_{k=1}^{k_n}Z_{n,k} )\big|\Big]
\le C\Capc\Big(\big|\sum_{k=1}^{k_n}Z_{n,k} \big|>N\Big)\\
\le  & C N^{-2}\sup_n \Sbep\Big[\big(\sum_{k=1}^{k_n}Z_{n,k}\big)^2\Big]
 \le CN^{-2} \sup_n \left(\Sbep\left[\delta_{k_n}\right]+\Sbep\left[\chi_{k_n}^2\right]\right) \\
 \le& 3C N^{-2}\to 0\;\; \text{ as } N\to \infty
\end{align*}
by (\ref{eqvariance}). The proof of Theorem \ref{thCLTM} is now completed. $\Box$

\subsection{Proof of the functional central limit theorem}

For proving the functional central limit theorem,
 we need a more lemma.

\begin{lemma}\label{lemfinite4}   Suppose that the operators $\Sbep_{n,k}$ satisfy (a) and (b), $\bm X_n\in \mathscr{H}_{n,k_n^{\prime}}\subset \mathscr{H}$ is a $d_1$-dimensional random vector, and    $\bm Y_n\in  \mathscr{H}$ is a $d_2$-dimensional random vector. Write $\mathscr{H}_{n}=\mathscr{H}_{n,k_n^{\prime}}$.
 Assume  that $\bm X_n\overset{d}\to \bm X$, and for any  bounded Lipschitz function $\varphi(\bm x,\bm y):\mathbb R_{d_1}\bigotimes \mathbb R_{d_2}\to \mathbb R$,
\begin{equation}\label{eqlemfinite4.1} \Sbep\left[\Big|\Sbep[\varphi(\bm x,\bm Y_n)|\mathscr{H}_n]-\widetilde{\mathbb E}[\varphi(\bm x,\bm Y)]\Big|\right]\to 0,\;\; \forall \bm x,
\end{equation}
where $\bm X$, $\bm Y$ are two random vectors in a sub-linear expectation space  $(\Omega, \mathscr{H}, \widetilde{\mathbb E})$ with $\widetilde{\Capc}(\|\bm X\|>\lambda)\to 0$  and $\widetilde{\Capc}(\|\bm Y\|>\lambda)\to 0$  as $\lambda\to \infty$. Then
\begin{equation}\label{eqlemfinite4.2} (\bm X_n,\bm Y_n)\overset{d }\to  (\widetilde{\bm X},\widetilde{\bm Y}),
\end{equation}
 where   $\widetilde{\bm Y}$ is independent to  $\widetilde{\bm X}$, $\widetilde{\bm X}\overset{d}= \bm X$ and $\widetilde{\bm Y}\overset{d}= \bm Y$.
\end{lemma}
{\bf Proof.} Suppose  $\varphi(\bm x,\bm y):\mathbb R_{d_1}\bigotimes \mathbb R_{d_2}\to \mathbb R$ is a bounded continuous function. We want to show that
\begin{equation}\label{eqlemfinite4.3} \Sbep\left[\varphi(\bm X_n,\bm Y_n)\right]\to \widetilde{\mathbb E}\left[\varphi(\widetilde{\bm X},\widetilde{\bm Y})\right].
\end{equation}
 First we assume that $\varphi(\bm x,\bm y)$ is a bounded Lipschitz function.    Without loss of generality, we assume $0\le \varphi(\bm x,\bm y)\le 1$ and $\left|\varphi(\bm x_1,\bm y_1)-\varphi(\bm x_2,\bm y_2)\right|\le \|\bm x_1-\bm x_2\|+\|\bm y_1-\bm y_2\|$. Let $ g_n(\bm x)=\Sbep\left[\varphi(\bm x, \bm Y_n)\big|\mathscr{H}_n\right]$
 and $ g(\bm x)=\widetilde{\mathbb E}[\varphi(\bm x,  \widetilde{\bm Y})]$. Then
 $$ |g(\bm x_1)-g(\bm x_2)|\le \widetilde{\mathbb E}\big[|\varphi(\bm x_1,  \widetilde{\bm Y})-\varphi(\bm x_2,  \widetilde{\bm Y})|\big]\le \|\bm x_1-\bm x_2\| $$
 and
 $$\left|\Sbep[\varphi(\bm X_n,\bm Y_n)|\mathscr{H}_n]-g_n(\bm x)\right|\le  \Sbep\big[|\varphi(\bm X_n,\bm Y_n)-\varphi(\bm x,\bm Y_n)|\big|\mathscr{H}_n\big]
 \le \|\bm X_n-\bm x\| \;\; \text{ in } L_1, $$
 by Lemma \ref{lemma4.0.2}.
 We use an argument of Hu, Li and Liu (2018) (c.f.  Proposition 3.4) to approximate  the  function $\varphi(\bm x,\bm y)$. For fixed $N\ge 1$, denote $B_N(0)=\{\bm x: \|x\|\le N\}$. By partition of
unity theorem, there exist $h_i\in C_{b,lip}(\mathbb R^{d_1})$, $i=1,\cdots, k_N$, such that $0\le h_i(\bm x)\le 1$, $I_{B_N(0)}\le \sum_{i=1}^{k_n}h_i(\bm x)\le 1$, and the diameter of support $\lambda(supp(h_i))\le 1/N$. Choose $\bm x_i$ such that $h_i(\bm x_i)>0$. Then
\begin{align*}
&\Big|\Sbep[\varphi(\bm X_n,\bm Y_n)|\mathscr{H}_n]-\sum_{i=1}^{k_N}h_i(\bm X_n)g_n(\bm x_i)\Big|\\
\le & \sum_{i=1}^{k_N}h_i(\bm X_n)\left|\Sbep[\varphi(\bm X_n,\bm Y_n)|\mathscr{H}_n]-g_n(x_i)\right|
  +\big(1- \sum_{i=1}^{k_N}h_i(\bm X_n)\big)\Big|\Sbep[\varphi(\bm X_n,\bm Y_n)|\mathscr{H}_n]\Big|\\
\le & \sum_{i=1}^{k_N}h_i(\bm X_n)\|\bm X_n-\bm x_i\|+\big(1- \sum_{i=1}^{k_N}h_i(\bm X_n)\big)
\le \frac{1}{N}+\big(1- \sum_{i=1}^{k_N}h_i(\bm X_n)\big) \;\; \text{ in } L_1.
\end{align*}
It follows that
\begin{align*}
&\Big|\Sbep[\varphi(\bm X_n,\bm Y_n)]-\Sbep\big[\sum_{i=1}^{k_N}h_i(\bm X_n)g_n(x_i)\big]\Big| \\
= &\Big|\Sbep\left[\Sbep[\varphi(\bm X_n,\bm Y_n)|\mathscr{H}_n]\right]-\Sbep\big[\sum_{i=1}^{k_N}h_i(\bm X_n)g_n(x_i)\big]\Big|\\
\le &\Sbep\left[\Big|\Sbep[\varphi(\bm X_n,\bm Y_n)|\mathscr{H}_n]-\sum_{i=1}^{k_N}h_i(\bm X_n)g_n(\bm x_i)\Big|\right] \\
\le & \frac{1}{N}+\Sbep\big[1- \sum_{i=1}^{k_N}h_i(\bm X_n)\big].
\end{align*}
Similarly,
\begin{align*}
&\Big|\widetilde{\mathbb E}\big[\varphi(\widetilde{\bm X},\widetilde{\bm Y})\big]-\widetilde{\mathbb E}\big[\sum_{i=1}^{k_N}h_i(\widetilde{\bm X})g(\bm x_i)\big]\Big|=\Big|\widetilde{\mathbb E}\big[g(\widetilde{\bm X})\big]-\widetilde{\mathbb E}\big[\sum_{i=1}^{k_N}h_i(\widetilde{\bm X})g(\bm x_i)\big]\Big| \\
\le & \widetilde{\mathbb E}\Big[ \big|g(\widetilde{\bm X})-\sum_{i=1}^{k_N}h_i(\widetilde{\bm X})g(\bm x_i)\big|\Big]
\le  \frac{1}{N}+\widetilde{\mathbb E}\big[1- \sum_{i=1}^{k_N}h_i(\widetilde{\bm X})\big].
\end{align*}
On the other hand, we have
\begin{align*}
& \left|\Sbep\big[\sum_{i=1}^{k_N}h_i(\bm X_n)g_n(\bm x_i)\big]-\Sbep\big[\sum_{i=1}^{k_N}h_i(\bm X_n)g(\bm x_i)\big]\right|\\
& \quad \le   \sum_{i=1}^{k_N}\Sbep[|g_n(\bm x_i)-g(\bm x_i)|] \text{ as } n\to \infty,
\end{align*}
 by (\ref{eqlemfinite4.1}), and
 $$  \Sbep\big[\sum_{i=1}^{k_N}h_i(\bm X_n)g(\bm x_i)\big] \to \widetilde{\mathbb E}\big[\sum_{i=1}^{k_N}h_i(\widetilde{\bm X})g(\bm x_i)\big], $$
 $$ \Sbep\big[1- \sum_{i=1}^{k_N}h_i(\bm X_n)\big]\to \widetilde{\mathbb E}\big[1- \sum_{i=1}^{k_N}h_i(\widetilde{\bm X})\big] $$
as $n\to \infty$, by the fact that $\bm X_n\overset{d}\to \widetilde{\bm X}$.  Combining the above arguments yields
\begin{align*}
&\limsup_{n\to \infty}\Big|\Sbep[\varphi(\bm X_n,\bm Y_n)]-\widetilde{\mathbb E}[\varphi(\widetilde{\bm X},\widetilde{\bm Y})]\Big|\\
\le &\frac{2}{N}+2\widetilde{\mathbb E}\big[1- \sum_{i=1}^{k_N}h_i(\widetilde{\bm X})\big]
\le   \frac{2}{N}+2\widetilde{\Capc}\left(\|\bm X\|>N\right)\to 0 \text{ as } N\to \infty.
\end{align*}
Hence (\ref{eqlemfinite4.3}) is proved for any bounded Lipschitz function $\varphi$.
For a bounded and uniformly continuous function $\varphi$, we define
$$\varphi_{\delta}=\varphi\ast \psi_{\delta} \text{ with } \psi_{\delta}(\bm x,\bm y)=\frac{1}{(2\pi \delta)^{(d_1+d_2)/2}}\exp\left\{-\frac{\sum_{i=1}^{d_1}x_i^2+\sum_{j=1}^{d_2}y_j^2}{2\delta}\right\}.
$$
Then $\varphi_{\delta}$ is a bounded Lipschitz function with $\sup_{\bm x,\bm y}|\varphi_{\delta}(\bm x,\bm y)-\varphi(\bm x,\bm y)|\to 0$ as $\delta\to 0$. Hence,
(\ref{eqlemfinite4.3}) holds for any  bounded and uniformly continuous function $\varphi$.
 Finally, let $\varphi(\bm x,\bm y)$ be a  bounded continuous function with $|\varphi(\bm x,\bm y)|\le M$. Let $\lambda>0$. For $\bm x=(x_1,\ldots, x_d)$, denote $\bm x_{\lambda}=\big((-\lambda)\vee(x_1\wedge \lambda)\lambda,\ldots, (-\lambda)\vee(x_{d_1}\wedge \lambda)\big)$ and define $\bm y_{\lambda}$ similarly. Let $\varphi_{\lambda}(\bm x,\bm y)=\varphi(\bm x_{\lambda},\bm y_{\lambda})$. Then $\varphi_{\lambda}$ is a bounded uniformly continuous function with
$$|\varphi_{\lambda}(\bm x,\bm y)-\varphi(\bm x,\bm y)|\le 2MI\{\|\bm x\|>\lambda\}+2MI\{\|\bm y\|>\lambda\}. $$
It follows that
\begin{align*}
& \limsup_{n\to \infty} \left|\Sbep[\varphi(\bm X_n,\bm Y_n)]-\widetilde{\mathbb E}[\varphi(\widetilde{\bm X},\widetilde{\bm Y})\right| \\
\le & \limsup_{n\to \infty} \left|\Sbep[\varphi_{\lambda}(\bm X_n,\bm Y_n)]-\widetilde{\mathbb E}[\varphi_{\lambda}(\widetilde{\bm X},\widetilde{\bm Y})]\right| \\
& + 2M\limsup_{n\to \infty}
\big\{\Capc(\|\bm X_n\|>\lambda)+\Capc(\|\bm Y_n\|>\lambda)\big\}\\
& +
2M
\big\{\widetilde{\Capc}(\|\bm X\|>\lambda)+\widetilde{\Capc}(\|\bm Y\|>\lambda)\big\}\\
\le & 4M\big\{\widetilde{\Capc}(\|\bm X\|>\lambda/2)+\widetilde{\Capc}(\|\bm Y\|>\lambda/2)\big\}\to 0\; \text{ as } \lambda\to \infty.
\end{align*}
 The proof is completed. $\Box$

 \begin{remark} In the original proofs of Lemma \ref{lemfinite4} and Theorem \ref{thFCLTM}, we need an additional assumption on the   operators $\Sbep_{n,k}$  as follows.
\begin{description}
\item[\rm (a$^{\prime}$)] If $\bm X=(X_1,\ldots, X_d)\in \mathscr{H}_{n,k}$, $Z\in \mathscr{H}$ and $\varphi(\bm x,y)$ is a bounded Lipschitz
function, then
$$\Sbep[\varphi(\bm X, Z)]=\Sbep\left[\Sbep_{n,k}\left[\varphi(\bm x,Z)\right]\Big|_{\bm x=\bm X}\right]. $$
 \end{description}
 \end{remark}
 We thank one of the referees mentioning  us the Proposition 3.4 of Hu, Li and Liu (2018) which helps us to remove this condition, though we fail to verify this proposition when the point by point monotonicity of the conditional sub-linear expectation (c.f. Definition 3.1 (1) of Hu, Li and Liu (2018)) is replaced by the $L_1$-monotonicity (c.f. Lemma 4.3 (c)).

 \bigskip
 {\bf Proof of Theorem \ref{thFCLTM}}. With the same argument as that at the beginning of the proof of Theorem \ref{thCLTM},   we can assume that $\delta_{k_n}=\sum_{k=1}^{k_n}\Sbep[Z_{n,k}^2|\mathscr{H}_{n,k-1}]\le 2\rho(1) $ in $L_1$, $\chi_{k_n}=:\sum_{k=1}^{k_n}\left\{ |\Sbep[Z_{n,k} |\mathscr{H}_{n,k-1}]|+|\cSbep[Z_{n,k} |\mathscr{H}_{n,k-1}]|\right\}<1$ in $L_1$ and $|Z_{n,k}|\le \epsilon_n$, $k=1,\ldots,k_n$, with a sequence $0<\epsilon_n\to 0$. Let $0<t_1<t_2\le 1$. Consider $\{Z_{n,k}^{\ast}=:Z_{n,\tau_n(t_1)+k}; k= 1,\ldots, k_n^{\ast}\}$, $S_i^{\ast}=\sum_{k=1}^iZ_{n,\tau_n(t_1)+k}$,
 $k_n^{\ast}=\tau_n(t_2) -\tau_n(t_1)$.
 Then $S_{k_n^{\ast}}^{\ast}=S_{n, \tau_n(t_2)}-S_{n, \tau_n(t_1)}=\sum_{k=1}^{k_n^{\ast}}Z_{n,\tau_n(t_1)+k}$,
 $$ \sum_{k=1}^{k_n^{\ast}}\Sbep\left[Z_{n,\tau_n(t_1)+k}^2\Big|\mathscr{H}_{n,\tau_n(t_1)+k-1}\right]\overset{\Capc}\to \rho (t_2)-\rho(t_1). $$
 By Theorem \ref{thCLT},
 $$ S_{n, \tau_n(t_2)}-S_{n, \tau_n(t_1)} \overset{d}\to W(\rho(t_2))-W(\rho(t_1)). $$
 Further, for any a bounded Lipschitz   function $\varphi(\bm u, x)$, let $V^{\bm u} (t, x)$ be the unique viscosity solution of the following equation,
$$ \partial_t V^{\bm u} + G( \partial_{xx}^2 V^{\bm u})=0,\;\;  (t, x) \in [0,\varrho+ h] \times \mathbb R, \; V^{\bm u}|_{t=\varrho+h} = \varphi(\bm u, x),
$$
where $\varrho=\rho(t_2)-\rho(t_1)$. With the same argument for showing (\ref{eqconvergence}), we can show that
\begin{equation}\label{eqproofFCLTM.2}
\Sbep\left[\left|\Sbep\left[\varphi\big(\bm u, S_{n, \tau_n(t_2)}-S_{n, \tau_n(t_1)}\big)\big|\mathscr{H}_{n,\tau_n(t_1)]}\right]
-\widetilde{\mathbb E}\left[\varphi\big(\bm u,  W(\rho(t_2))-W(\rho(t_1)) \big)\right]\right|\right]\to 0.
\end{equation}
The only difference is that (\ref{eqtheorem2proof.2}), (\ref{eqtheorem2proof.3}) and (\ref{eqtheorem2proof.4}) are needed to be replaced, respectively, by
\begin{align*}
 &\Sbep\left[\Big|\Sbep \Big[V^{\bm u}(\delta_{k_n^{\ast}}^{\ast}\wedge h_0, S_{k_n^{\ast}}^{\ast})\Big|\mathscr{H}_{n,\tau_n(t_1)}\Big]-V^{\bm u}(0,0)-\Sbep\Big[\sum_{i=0}^{k_n^{\ast}-1}(J_{n,1,\ast}^i+J_{n,2,\ast}^i)\Big|\mathscr{H}_{n,\tau_n(t_1)}\Big]\Big|\right] \\
 &\quad  \le
 \Sbep\Big[\big|V^{\bm u}(\delta_{k_n^{\ast}}^{\ast}\wedge h_0, S_{k_n^{\ast}}^{\ast}) -V^{\bm u}(0,0)-\sum_{i=0}^{k_n^{\ast}-1}(J_{n,1,\ast}^i+J_{n,2,\ast}^i)\big|\Big],
\end{align*}
\begin{align*}   \Sbep\Big[\sum_{i=0}^{k_n^{\ast}-1}J_{n,1,\ast}^i\Big|\mathscr{H}_{n,\tau_n(t_1)}\Big]
=&\Sbep\Big[\Sbep\Big[\sum_{i=0}^{k_n^{\ast}-1}J_{n,1,\ast}^i\Big|\mathscr{H}_{n,\tau_n(t_1)+k_n^{\ast}-1}\Big]\Big|\mathscr{H}_{n,\tau_n(t_1)}\Big]
\\
=& \Sbep\Big[\sum_{i=0}^{k_n^{\ast}-2}J_{n,1,\ast}^i+\Sbep\Big[J_{n,1,\ast}^{k_n-1}\big|\mathscr{H}_{n,\tau_n(t_1)+k_n^{\ast}-1}\Big]\Big|\mathscr{H}_{n,\tau_n(t_1)}\Big]
\\
=&
\Sbep\Big[\sum_{i=0}^{k_n^{\ast}-2}J_{n,1,\ast}^i\Big|\mathscr{H}_{n,\tau_n(t_1)}\Big]=\ldots=0 \; \text{ in } L_1,
\end{align*}
and
 \begin{align*}  &\Sbep\Big[\sum_{i=0}^{k_n-1}(\pm J_{n,2,\ast}^i-\widetilde{J}_{n,2,\ast}^i)\Big|\mathscr{H}_{n,\tau_n(t_1)}\Big] \\
=&\Sbep\left[\Sbep\Big[\sum_{i=0}^{k_n-2}(\pm J_{n,2,\ast}^i-\widetilde{J}_{n,2,\ast}^i)+\Sbep\left[\pm J_{n,1}^{k_n-1}-\widetilde{J}_{n,2}^{k_n-1}\big|\mathscr{H}_{n,k_n-1}\right]\Big|\mathscr{H}_{n,\tau_n(t_1)}\Big]\right]\\
\le & \Sbep\Big[\sum_{i=0}^{k_n-2}(\pm J_{n,2,\ast}^i-\widetilde{J}_{n,2,\ast}^i)\Big|\mathscr{H}_{n,\tau_n(t_1)}\Big]\le \ldots\le 0\; \text{ in } L_1.
\end{align*}
where $J_{n,1,\ast}^i$, $J_{n,2,\ast}^i$ and $\widetilde{J}_{n,2,\ast}^i$ are defined the same as $J_{n,1}^i$, $J_{n,2}^i$ and $\widetilde{J}_{n,2}^i$ with $\{Z_{n,k}^{\ast}\}$ taking the place of $\{Z_{n,k}\}$.
On the other hand, note
$ S_{n,\tau_n(t_1)}\overset{d}\to W(\rho(t_1)). $
Hence,
$$ \Big(S_{n,\tau_n(t_1)}, S_{n,\tau_n(t_2)}-S_{n,\tau_n(t_1)}\Big)\overset{d}\to \Big(W(\rho(t_1)), W(\rho(t_2))-W(\rho(t_1))\Big). $$
by (\ref{eqproofFCLTM.2}) and Lemma  \ref{lemfinite4}. By induction, for any $0=t_0<\ldots< t_d\le 1$,
\begin{align*}
&  \Big( S_{n,\tau_n(t_1)}-S_{n,\tau_n(t_0)},\cdots, S_{n,\tau_n(t_d)}-S_{n,\tau_n(t_{d-1})}\Big)\\
  \overset{d}
\to  & \Big( W(\rho(t_1))-W(\rho(t_0)),\cdots, W(\rho(t_d))-W(\rho(t_{d-1}))\Big),
\end{align*}
which implies (\ref{eqfinitedimension}).
So, we have shown the convergence of finite dimensional distributions of $W_n$. By Theorem 9 of Peng (2010) on the tightness  and the argument of   Lin and Zhang (2017) or Zhang (2015), to show that (\ref{eqFCLTM})   holds for bounded   continuous function $\varphi$, it is sufficient to show that for any $\epsilon^{\prime}>0$,
\begin{equation}\label{eqpoofthFCLTMtight} \lim_{\delta\to 0} \limsup_{n\to \infty} \Capc\left( w_{\delta}\left(W_n\right)\ge 3\epsilon^{\prime}\right)=0,
\end{equation}
where $\omega_{\delta}(x)=\sup_{|t-s|<\delta,t,s\in[0,1]}|x(t)-x(s)|$ (c.f. Proposition \ref{proptightness}  in Appendix \ref{appendixB}).
 Assume $0<\delta<1/10$. Let $0=t_0<t_1\ldots<t_K=1$ such that $t_k-t_{k-1}=\delta$, and let $t_{K+1}=t_{K+2}=1$. It is easily seen that
$$ \Capc\left( w_{\delta}\left(W_n\right)\ge 3\epsilon^{\prime}\right)
\le 2\sum_{k=0}^{K-1}   \Capc\left( \max_{s\in [t_k,t_{k+2}]} |S_{n,\tau_n(s)}-S_{n,\tau_n(t_k)}|\ge  \epsilon^{\prime}\right). $$
On the other hand,  for $t,\gamma>0$,    by (\ref{eqlemRosenIeq.4}) we have
 \begin{align*}
 & \Sbep \left[ \max_{s\le \gamma} |S_{n,\tau_n(t+s)}-S_{n,\tau_n(t)}|^4\right]\nonumber\\
 \le &  C \Sbep\left[\sum_{k=\tau_n(t)+1}^{\tau_n(t+\gamma)}\Sbep\left[Z_{n,k}^4\big|\mathscr{H}_{n,k-1}\right]\right]
  +C  \Sbep\left[\Big(\sum_{k=\tau_n(t)+1}^{\tau_n(t+\gamma)}\Sbep\left[Z_{n,k}^2\big|\mathscr{H}_{n,k-1}\right]\Big)^2\right]\nonumber \\
 &+C\Sbep\left[\Big(\sum_{k=\tau_n(t)+1}^{\tau_n(t+\gamma)}\left\{ |\Sbep[Z_{n,k} |\mathscr{H}_{n,k-1}]|+|\cSbep[Z_{n,k} |\mathscr{H}_{n,k-1}]|\right\}\Big)^4\right]\nonumber\\
 \le &C  \Sbep\left[\Big(\sum_{k=\tau_n(t)+1}^{\tau_n(t+\gamma)}\Sbep\left[Z_{n,k}^2\big|\mathscr{H}_{n,k-1}\right]\Big)^2\right]
 +C\epsilon_n^2 \cdot 2\rho +C\Sbep[\chi_{k_n}^4].
 \end{align*}
 The last two terms above will go to zero by (\ref{eqCLTCondM.4}). For considering the first   term, we note
\begin{equation}\label{eqproofFCLTM.8}
2\rho(1)\ge   \sum_{k=\tau_n(t)+1}^{\tau_n(t+\gamma)}\Sbep\left[Z_{n,k}^2\big|\mathscr{H}_{n,k-1}\right] \overset{\Capc}\to   \rho(t+\gamma)-\rho(t).
\end{equation}
 It follows that
$$\Sbep\left[ \Big(\sum_{k=\tau_n(t)+1}^{\tau_n(t+\gamma)}\Sbep\left[Z_{n,k}^2\big|\mathscr{H}_{n,k-1}\right]\Big)^2\right]\to (\rho(t+\gamma)-\rho(t))^2. $$
 So, we conclude that
 \begin{align*}
& \limsup_n   2\sum_{k=0}^{K-1}   \Capc\left( \max_{s\in [t_k,t_{k+2}]} |S_{n,\tau_n(s)}-S_{n,\tau_n(t_k)}|\ge  \epsilon^{\prime}\right)\\
\le &   \limsup_n   2\sum_{k=0}^{K-1} \Big(\frac{1}{\epsilon ^{\ast}}\Big)^4  \Sbep\left[ \max_{s\in [t_k,t_{k+2}]} |S_{n,\tau_n(s)}-S_{n,\tau_n(t_k)}|^4 \right]\\
\le & C\sum_{k=0}^{K-1} \frac{1}{ (\epsilon ^{\ast})^4}  \big(\rho(t_{k+2})-\rho(t_k)\big)^2
\le   C  \frac{\rho(1)}{ (\epsilon ^{\ast})^4}   \sup_{|t-s|\le 2\delta} \big|\rho(t)-\rho(s)\big| \to 0
 \end{align*}
 by taking $\delta\to 0$.  Hence, (\ref{eqpoofthFCLTMtight}) is verified. And the proof is completed. $\Box$


\begin{appendix}

\section*{Appendix}

\section{The properties of the conditional expectations. }\label{appendixA}
\setcounter{equation}{0}

In this appendix, we give the proofs of Lemmas \ref{lemma4.0.1} and \ref{lemma4.0.2} on the  properties of the conditional expectation.

{\bf Proof of Lemma \ref{lemma4.0.1}.}
(1) is obvious. For (2), note that
 $$ \Sbep[((X-Y)^+)^q]\le \epsilon^q+c^q\Capc(X-Y\ge \epsilon)+ \Sbep\left[((X-Y-c)^+)^q\right] $$
 and
 $$\Sbep\left[((X-Y-c)^+)^q\right] \le \frac{\Sbep[((X-Y)^+)^p]}{c^{p-q}}\to 0 \text{ as } c\to \infty. $$
 The results follows.

 For (3), let $\epsilon>0$ and $M>0$ be given. Let $0<\delta<1$ such that $|x-y|\le \delta$ and $|y|\le M$ implies $|f(x)-f(y)|\le \epsilon$. Then,
$$
 \Capc\left(f(X)-f(Y)\ge \epsilon\right)\le  \Capc\left(X-Y\ge \delta \right)+ \Capc\left|Y|\ge M\right).
$$
  The result follows.

  For (4), note for $y, x\ge 0$,  $x^p-y^p\le p x^{p-1}(x-y)$. So,
  $$ \Sbep[X^p]-\Sbep[Y^p]\le p\Sbep[X^{p-1}(X-Y)^+]\le p(\Sbep[X^p])^{1/q} \left(\Sbep\left[\big((X-Y)^+\big)^p\right]\right)^{1/p}=0.  $$

  For (5), note that the countable additivity of $\Sbep$ implies
  $$ \Sbep[((X-Y)^+)^p]\le \int_0^{\infty} \Capc\left(((X-Y)^+)^p>y\right)dy=\int_0^{\infty} \Capc\left( X-Y >y^{1/p}\right)dy $$
 (cf. Lemma 3.9 of Zhang (2016)).   The result follows.  $\Box$

  \bigskip

{\bf Proof of Lemma Lemma \ref{lemma4.0.2}.} (c) Let $0\le f\in \mathscr{H}_{n,k}$ be a bounded random variable. Then
\begin{align*}
&\Sbep\left[f\big(\Sbep_{n,k}[X]-   \Sbep_{n,k}[Y]\big)\right]=\Sbep\left[ \Sbep_{n,k}\big[fX -   \Sbep_{n,k}[fY]\big]\right]\\
=&   \Sbep \big[fX -   \Sbep_{n,k}[fY]\big]\le \Sbep \big[fX -fY+fY-   \Sbep_{n,k}[fY]\big] \\
\le &\Sbep [f(X -Y)^+]+\Sbep[fY-   \Sbep_{n,k}[fY]\big]\le \Sbep[fY-   \Sbep_{n,k}[fY]\big]\\
=& \Sbep\left[ \Sbep_{n,k}\big[fY -   \Sbep_{n,k}[fY]\big]\right]=\Sbep\left[ \Sbep_{n,k}[fY] -   \Sbep_{n,k}[fY] \right]=0,
\end{align*}
which will implies $\Sbep\left[\big(\Sbep_{n,k}[X]-   \Sbep_{n,k}[Y]\big)^+\right]=0$. In fact, let $Z=\Sbep_{n,k}[X]-   \Sbep_{n,k}[Y]$ and choose $f$ to be a bounded Lipschitz function of $Z$ such that $I\{Z\ge 2\epsilon\}\le f\le I\{Z\ge \epsilon\}$. Then,
$$ \Sbep[Z^+]\le 2\epsilon +\Sbep[fZ]\le 2\epsilon. $$

 (d) The second inequality is due to (c). For the first one, let $Z= \Sbep_{n,k}[X]-\Sbep_{n,k}[Y]- \Sbep_{n,k}[X-Y]$. With the same argument as in (c), it is sufficient to show that $\Sbep[fZ]\le 0$ for any bounded $0\le f\in \mathscr{H}_{n,k}$. Now,
 \begin{align*}
 \Sbep[fZ]=&\Sbep\left[\Sbep_{n,k}\Big[fX -\Sbep_{n,k}[fY]- \Sbep_{n,k}[fX-fY]\Big]\right]\\
 =& \Sbep \Big[fX -\Sbep_{n,k}[fY]- \Sbep_{n,k}[fX-fY]\Big]\\
  =&\Sbep \Big[\big(fY -\Sbep_{n,k}[fY]\big)+\big(fX-fY- \Sbep_{n,k}[fX-fY]\big)\Big] \\
  \le &\Sbep  \big[fY -\Sbep_{n,k}[fY]\big]+\Sbep[fX-fY- \Sbep_{n,k}[fX-fY]\big] =0.
 \end{align*}

 (e) Suppose $k<l$. Let $Z=\Sbep_{n,k}\left[\left[\Sbep_{n,l} [ X]\right]\right]-\Sbep_{n,  k} [ X]$ and $f\ge 0$ be a bounded random variable  in $\mathscr{H}_{n,k}$. Then,
 \begin{align*}
 \Sbep[fZ]=&\Sbep\left[\Sbep_{n,k}\Big[\big[\Sbep_{n,l} [fX]\big]\Big]-\Sbep_{n,  k} [f X]\right]\\
 =&\Sbep\left[\Sbep_{n,k}\Big[\Sbep_{n,l} \big[fX-\Sbep_{n,  k} [f X] \big]\Big]\right]= \Sbep \Big[ fX-\Sbep_{n,k} [fX]\Big] =0.
 \end{align*}
which will imply $\Sbep[Z^+]=0$. On the other hand, note $-Z\le \Sbep_{n,k}\Big[X-\Sbep_{n,l}[X]\Big]$ by Property (d). We have
$$ \Sbep[f(-Z)]\le \Sbep\left[\Sbep_{n,k}\Big[fX-\Sbep_{n,l}[fX]\Big]\right]=\Sbep\Big[fX-\Sbep_{n,l}[fX]\Big]=0, $$
which implies $\Sbep[(-Z)^+]=0$. So, $\Sbep[|Z|]=0$,

 (f) Let $Z= \Sbep_{n,k}[X]$ and $0\le f\in \mathscr{H}_{n,k}$ be a bounded random variable with $fZ^+=0$ and $|f|\le 1$. Then $Z\in \mathscr{L}(\mathscr{H})$. We first show that $f|Z|^p\in \mathscr{L}(\mathscr{H})$ for any $p\ge 1$. It is obvious that $f|Z|\in \mathscr{L}(\mathscr{H})$.  Assume that $k\ge 1$  is  an   integer, and $f|Z|^p\in \mathscr{L}(\mathscr{H})$ for $p\le k$. Let $p^{\prime}\ge k$, $p^{\prime}\le p<p^{\prime}+1$. Note
  $$0\le  f |X| |Z|^{p-1}\le \frac{p^{\prime}+1-p}{p^{\prime}}|X|^{\frac{p^{\prime}}{p^{\prime}+1-p}}+\frac{p-1}{p^{\prime}}f |Z|^{p^{\prime}}. $$
  Choosing $p^{\prime}=k$ yields $ X, f |Z|^{p-1}, f |X| |Z|^{p-1}\in \mathscr{L}(\mathscr{H})$. So by the Properties (a), (b) and (d),
  \begin{align*}
  \Sbep[f|Z|^p]=&\Sbep\left[f|Z|^{p-1}(-\Sbep_{n,k}[X])\right]=\Sbep\left[-\Sbep_{n,k}[Xf|Z|^{p-1}]\right]\\
  \le &\Sbep\left[\Sbep_{n,k}[|X|\cdot f|Z|^{p-1}]\right]=\Sbep[f|X||Z|^{p-1}]<\infty
  \end{align*}
    if $k\le p<k+1$. Choosing $p^{\prime}=k+1/2$ and repeating the same argument yield $\Sbep[f|Z|^p]<\infty$ if $k+1/2\le p<k+3/2$.
   So, $f|Z|^p\in \mathscr{L}(\mathscr{H})$ for $k\le p\le k+1$. By the induction, for any $p\ge 1$, $\Sbep[f|Z|^p]<\infty$
 which will imply $\Sbep[(Z^-)^p]<\infty$. And similarly by choosing $f$ such that $fZ^-=0$, we will have $\Sbep[(Z^+)^p]<\infty$. So, we have $\Sbep[|Z|^p]<\infty$ for any $p\ge 1$. Finally, by (c), $|Z|\le M$ in $L_1$. Hence, by Lemma \ref{lemma4.0.1} (2), the result follows. The proof is now completed.   $\Box$

\section{Tightness. }\label{appendixB}
\setcounter{equation}{0}

\begin{proposition}\label{proptightness}
Let $\{Z_{n,k}; k=1,\ldots, k_n\}$ be an array of random variables with $\Sbep[|Z_{n,k}|]<\infty$, $k=1,\ldots, k_n$, and $\tau_n(t)$ be a non-decreasing function in $D_{[0,1]}$ which takes integer values with $\tau_n(0)=0$, $\tau_n(1)=k_n$. Define  $S_{n,i}=\sum_{k=1}^i Z_{n,k}$,
\begin{equation}\label{eqtight.1} W_n(t)= S_{n, \tau_n(t)}.
\end{equation}
 Assume that for any $\epsilon>0$,
\begin{equation}\label{eqtight.2} \lim_{\delta\to 0} \limsup_{n\to \infty} \Capc\left( w_{\delta}\left(W_n\right)\ge \epsilon\right)=0,
\end{equation}
where $\omega_{\delta}(x)=\sup_{|t-s|<\delta,t,s\in[0,1]}|x(t)-x(s)|$. Then $\{W_n\}$ is tight in $D_{[0,1]}$ endowed the Skorohod topology, i.e., for any $\eta>0$, there exists a compact set $K$ in $D_{[0,1]}$   such that
\begin{equation}\label{eqtight.3} \sup_n\Capc\left(W_n\not\in K\right) \le \eta.
\end{equation}
Further, if (\ref{eqfinitedimension}) holds for any $0<t_1<\cdots,t_d\le 1$, then (\ref{eqFCLTM}) holds.
\end{proposition}

{\bf Proof.} The proof of the tightness is similar to that of the tightness of probability measures (c.f. Billingsley (1968)). The only difference we shall note is that $\Capc$ may be not countably additive and may be not continuous. For $T_0\subset [0,1]$, define
$$w(x,T_0)=\sup_{t,s\in T_0}|x(t)-x(s)|, $$
and
$$ w^{\prime}_{\delta}(x)=\inf_{t_i}\max_{1\le i\le \nu} w\big(x,[t_{i-1},t_i)\big),  $$
where the infimum extends over all sets $\{t_i\}$ with
$$ 0=t_0<t_1<\cdots<t_{\nu-1}<t_{\nu}=1, \;\; \min_{1\le i\le \nu}(t_i-t_{i-1}) >\delta. $$
Note $w_{\delta}^{\prime}(x)\le w_{2\delta}(x)$,
$$ |x(t)|\le |x(0)|+\sum_{i=1}^k |x(it/k)-x((i-1)t/k)|\le |x(0)|+k w_{1/k}(x), $$
and $W_n(0)=0$.  From (\ref{eqtight.2}) it follows that
\begin{equation}\label{eqprooftight2}
\lim_{a\to \infty}\limsup_{n\to \infty}\Capc\left(\sup_{t}|W_n(t)|>a\right)=0
\end{equation}
and
\begin{equation}\label{eqprooftight3} \lim_{\delta\to 0} \limsup_{n\to \infty} \Capc\left( w_{\delta}^{\prime}\left(W_n\right)\ge \epsilon\right)=0, \;\; \forall \epsilon>0.
\end{equation}
For fixed $n$, let $0<t_1^n<\cdots<t_{\nu-1}^n\le 1$ be the jump times of the step function $\tau_n(t)$, $t_0^n=0$, $t_{\nu}^n=1$. Then
$$ w\big(W_n,[t_{i-1}^n,t_i^n)\big)=0, \;\; i=1,\cdots,\nu.
 $$
 Let $\delta_0^n =\min_{1\le i\le \nu-1}(t_i^n-t_{i-1}^n)$ if $ t_{\nu-1}^n=1$, and $=\min_{1\le i\le \nu}(t_i^n-t_{i-1}^n)$ if $ t_{\nu-1}^n<1$. Then
 \begin{equation} \label{eqprooftight4} w_{\delta}^{\prime}(W_n)=0 \;\; \text{when}\;\; \delta<\delta_0^n.
 \end{equation}
 On the other hand, it is obvious that
 $$\lim_{a\to \infty} \Capc\left(\sup_{t}|W_n(t)|>a\right)\le \lim_{a\to \infty} \frac{\sum_{k=1}^{k_n}\Sbep[|Z_{n,k}|]}{a}=0. $$
Hence, (\ref{eqprooftight2}) and (\ref{eqprooftight3}) imply that
\begin{equation}\label{eqproof5}
\lim_{a\to \infty}\sup_n\Capc\left(\sup_{t}|W_n(t)|>a\right)=0
\end{equation}
and
\begin{equation}\label{eqprooftight6} \lim_{\delta\to 0} \sup_n\Capc\left( w_{\delta}^{\prime}\left(W_n\right)\ge \epsilon\right)=0, \;\; \forall \epsilon>0.
\end{equation}
Now, for any $\eta>0$ and a sequence $0<\epsilon_k\to 0$,  choose $a>0$ and $0<\delta_k\to 0$ such that
 \begin{align*}
  \sup_n\Capc\left(\sup_{t}|W_n(t)|>a\right)<\eta/2\;\;\text{ and }\;\;
  \sup_n\Capc\left( w_{\delta_k}^{\prime}\left(W_n\right)> \epsilon_k\right)<\eta/2^{k+1}.
\end{align*}
Now, let $B_0=\{x\in D_{[0,1]}: \sup_t|x(t)|\le a\}$,  $B_k=\{x\in D_{[0,1]}: w_{\delta_k}^{\prime}\left(x\right)\le \epsilon_k \}$ and $A=\bigcap_{k=0}^{\infty} B_k$.
Then $\sup_{x\in A}\sup_t|x(t)|\le a$ and $\lim_{\delta\to 0}\sup_{x\in A}w_{\delta}^{\prime}\left(x\right)=0$. By the Arzal\'a-Ascoli thorem, the closure of A is a compact set in $D_{[0,1]}$.
On the other hand, by noting (\ref{eqprooftight4}),
\begin{align*}
\{W_n\not\in  cl(A)\}\subset  & \big\{\sup_{t}|W_n(t)|>a\big\} \bigcup_{k=1}^{\infty} \big\{ w_{\delta_k}^{\prime}\left(W_n\right)> \epsilon_k\big\}\\
\subset  & \big\{\sup_{t}|W_n(t)|>a\big\} \bigcup_{k: \delta_k\ge \delta_0^n} \big\{ w_{\delta_k}^{\prime}\left(W_n\right)> \epsilon_k\big\}.
\end{align*}
By the (finite) sub-additivity of $\Capc$, it follows that
\begin{align} \label{eqproofthight7}\Capc\Big(W_n\not\in  cl(A)\Big)\le  & \Capc\left(\sup_{t}|W_n(t)|>a\right)+\sum_{k:\delta_k\ge \delta_0^n}\Capc\left( w_{\delta_k}^{\prime}\left(W_n\right)> \epsilon_k\right) \\
< & \eta/2+\sum_{k=1}^{\infty} \eta/2^{k+1}=\eta.\nonumber
\end{align}
The proof  of the tightness (\ref{eqtight.3}) is completed.

Now, consider the G-Brownian motion $W$. In Zhang (2015), it is proved that
$$ \lim_{\delta\to 0}   \widetilde{\Capc}\left( w_{\delta}\left(W\right)\ge \epsilon\right)=0\;\; \text{ for any } \epsilon>0.
$$
Note that $\rho(\cdot)$ is a uniformly continuous function on $[0,1]$. It follows that
$$ \lim_{\delta\to 0}   \widetilde{\Capc}\left( w_{\delta}\left(W\circ\rho\right)\ge \epsilon\right)=0\;\; \text{ for any } \epsilon>0.
$$
With the same argument as (\ref{eqproofthight7}) one can show that for any $\eta>0$, there exists a compact set $K$ in $D_{[0,1]}$ such that
$$ \widetilde{\Capc}\left(W\circ\rho\not\in K\right)<\eta. $$

For $0=t_0< t_1<t_2\ldots<t_{d-1}<t_d= 1$, we define the projection $\pi_{t_1,\ldots,t_d}$ from $D_{[0,1]}$ to $\mathbb R^d$ by
$$ \pi_{t_1,\ldots,t_d}x=(x(t_1),\ldots, x(t_d)), $$
and define a map $\Pi^{-1}_{t_1,\ldots,t_d}$ from $\mathbb R^d$ to $D_{[0,1]}$  by
$$\Pi^{-1}_{t_1,\ldots,t_d}(x_1,\ldots,x_d)= \begin{cases}
0,  \; \text{ if }  t\in \big[t_0,t_1\big);\;\;
 x_k,  \; \text{ if } t\in \big[t_k, t_{k+1}\big) \; (k=1,\ldots, d);\\
 x_d,  \; \text{ if } t=t_d.
\end{cases}$$
Then  $\Pi_{t_1,\ldots,t_d}^{-1}$ is a continuous map. Denote $\widetilde{\pi}_{t_1,\ldots,t_d}=\Pi^{-1}_{t_1,\ldots,t_d}\circ \pi_{t_1,\ldots,t_d}$.
Let $\varphi\in C_b\big(D_{[0,1]})$. Then $\varphi(\widetilde{\pi}_{t_1,\ldots,t_d} x)=\varphi\circ\Pi^{-1}_{t_1,\ldots,t_d}(x(t_1),\ldots, x(t_d))$ and $\varphi\circ\Pi^{-1}_{t_1,\ldots,t_d}\in C_b(\mathbb R^d)$. By  (\ref{eqfinitedimension}) on the convergence of the finite-dimensional distributions of $W_n$, it follows that
\begin{align*}
&\lim_{n\to \infty}\Sbep\left[\varphi\left(\widetilde{\pi}_{t_1,\ldots,t_d} W_n\right)\right]
=\lim_{n\to \infty} \Sbep\left[\varphi\circ\Pi^{-1}_{t_1,\ldots,t_d}\left(W_n(t_1),\ldots,W_n(t_d)\right)\right]\\
& \qquad =\widetilde{\mathbb E}\left[\varphi\circ\Pi^{-1}_{t_1,\ldots,t_d}\left(W(\rho(t_1)),\ldots,W(\rho(t_d))\right)\right]
=\widetilde{\mathbb E}\left[\varphi\left(\widetilde{\pi}_{t_1,\ldots,t_d} W\circ\rho\right)\right].
\end{align*}
Now,  suppose that $t_{i+1}-t_i<\delta$ for $i=0,\ldots, d-1$. Recall  $\omega_{\delta}(x)=\sup\limits_{|t-s|<\delta}|x(t)-x(s)|$, and let $d_0(\cdot,\cdot)$ be the  Skorohod distance in $D_{[0,1]}$ and $\|x\|=\sup\limits_{0\le t\le 1}|x(t)|$. It is easily seen that
$  d_0\left(\widetilde{\pi}_{t_1,\ldots,t_d}x, x\right)\le \left\|\widetilde{\pi}_{t_1,\ldots,t_d}x -x \right\|\le \omega_{\delta}(x). $
Let $\epsilon>0$ be given. Since  $\varphi$ is a continuous function, for each $x$, there is an $\epsilon_x>0$ such that
$$ \left|\varphi(x)-\varphi(y)\right|<\epsilon/2 \text{ whenever } d_0(x,y)<\epsilon_x. $$
Let $K\subset D_{[0,1]}$ be a compact set. Then it can be covered by a union of finite many of the sets $\{y: d_0(x,y)<\epsilon_x/2\}$, $x\in K$. So, there is an $\epsilon_K>0$ such that
 $ \left|\varphi(x)-\varphi(y)\right|<\epsilon$  whenever $d_0(x,y)<\epsilon_K$ and $x\in K$.
Denote $M=\sup_x|\varphi(x)|$. It follows that
 \begin{align*}
 \left|\varphi\left(\widetilde{\pi}_{t_1,\ldots,t_d} x\right)-\varphi(x)\right|
 <\epsilon+2M I\{\omega_{\delta}(x)\ge \epsilon_K\}+2MI\{x\not\in K\}.
 \end{align*}
  By   the tightness of $\{W_n\}$ and $W\circ\rho$, respectively,  we can choose $K$ and $\delta$ such that
 $$ \sup_n\Capc\left(\omega_{\delta}(W_n)\ge \epsilon_K\right)+  \sup_n\Capc\left(W_n\not\in K\right)
 \le \frac{\epsilon}{4M} \;\; \text{and }$$
 $$\widetilde{\Capc}\left(\omega_{\delta}(W\circ\rho)\ge \epsilon_K\right)+  \widetilde{\Capc}\left(W\circ\rho\not\in K\right)
 \le \frac{\epsilon}{4M}. $$
 Hence
 \begin{align*}
 & \left| \Sbep\left[ \varphi(W_n)\right]- \widetilde{\mathbb E}\left[ \varphi(W\circ\rho)\right]\right|\\
 \le & \left| \Sbep\left[ \varphi\big(\widetilde{\pi}_{t_1,\ldots, t_d}W_n\big)\right]- \widetilde{\mathbb E}\left[\varphi\big(\widetilde{\pi}_{t_1,\ldots, t_d}W\circ\rho\big)\right]\right| \\
 &+\left|\Sbep\left[\varphi(W_n)\right]-\Sbep\left[\varphi\left(\widetilde{\pi}_{t_1,\ldots,t_d} W_n\right)\right]\right|
 +\left|\widetilde{\mathbb E}\left[\varphi\big(\widetilde{\pi}_{t_1,\ldots, t_d}W\circ\rho\big)\right]-\widetilde{\mathbb E}\left[ \varphi(W\circ\rho)\right]\right|\\
 \le & \left| \Sbep\left[ \varphi\big(\widetilde{\pi}_{t_1,\ldots, t_d}W_n\big)\right]- \widetilde{\mathbb E}\left[\varphi\big(\widetilde{\pi}_{t_1,\ldots, t_d}W\circ\rho\big)\right]\right| \\
 &+2\epsilon+ 2M \Capc\left(\omega_{\delta}(W_n)\ge \epsilon_K\right)+2M \Capc\left(W_n\not\in K\right) \\
 &+ 2M \widetilde{\Capc}\left(\omega_{\delta}(W\circ\rho)\ge \epsilon_K\right)+2M \widetilde{\Capc}\left(W\circ\rho\not\in K\right)\\
 \le & \left| \Sbep\left[ \varphi\big(\widetilde{\pi}_{t_1,\ldots, t_d}W_n\big)\right]- \widetilde{\mathbb E}\left[\varphi\big(\widetilde{\pi}_{t_1,\ldots, t_d}W\circ\rho\big)\right]\right|+3\epsilon.
 \end{align*}
Letting $n\to \infty$ and then $\epsilon\to 0$ completes the proof of (\ref{eqFCLTM}).
$\Box$

\end{appendix}
\section*{Acknowledgements} Special thanks go to the anonymous referees and the associate editor for    their
constructive comments, which led to a much improved version of this paper.


\end{document}